\documentclass[11pt,reqno]{amsart}
\usepackage{enumerate}
\usepackage{amsfonts}
\usepackage{amsmath}
\usepackage{amsthm}
\usepackage{amssymb}
\usepackage{latexsym}
\usepackage{multicol}
\usepackage{verbatim}
\usepackage{amscd}
\usepackage{hyperref}
\usepackage{pb-diagram}

\advance\textwidth by 1.2in \advance\oddsidemargin by -.6in \advance\evensidemargin by -.6in
\newtheorem*{cor}{Corollary}%[section]
\newtheorem*{lem}{Lemma}
\newtheorem*{prop}{Proposition}

\theoremstyle{definition}
\newtheorem*{defn}{Definition}
\theoremstyle{definition}
\newtheorem{thm}{Theorem}

\newtheorem*{rem}{Remark}

\newenvironment{pf}{\proof}{\endproof}
\newcounter{cnt}
 \makeatletter
\def\mydggeometry{\makeatletter\dg@YGRID=1\dg@XGRID=20\unitlength=0.003pt\makeatother}
\makeatother \theoremstyle{remark}

\numberwithin{equation}{section}

\begin{document}

\newcommand{\thmref}[1]{Theorem~\ref{#1}}
\newcommand{\secref}[1]{Section~\ref{#1}}
\newcommand{\lemref}[1]{Lemma~\ref{#1}}
\newcommand{\propref}[1]{Proposition~\ref{#1}}
\newcommand{\corref}[1]{Corollary~\ref{#1}}
\newcommand{\remref}[1]{Remark~\ref{#1}}
\newcommand{\defref}[1]{Definition~\ref{#1}}
\newcommand{\er}[1]{(\ref{#1})}
\newcommand{\id}{\operatorname{id}}
\newcommand{\wt}{\operatorname{wt}}
\newcommand{\tensor}{\otimes}
\newcommand{\nc}{\newcommand}
\newcommand{\rnc}{\renewcommand}
\newcommand{\qbinom}[2]{\genfrac[]{0pt}0{#1}{#2}}
\nc{\cal}{\mathcal} \nc{\goth}{\mathfrak} \rnc{\bold}{\mathbf}
\renewcommand{\frak}{\mathfrak}
\newcommand{\supp}{\operatorname{supp}}
\newcommand{\ad}{\operatorname{ad}}
\newcommand{\Maj}{\operatorname{Maj}}
\renewcommand{\Bbb}{\mathbb}
\nc\bomega{{\mbox{\boldmath $\omega$}}} \nc\bpsi{{\mbox{\boldmath $\Psi$}}}
 \nc\balpha{{\mbox{\boldmath $\alpha$}}}
 \nc\bpi{{\mbox{\boldmath $\pi$}}}

\newcommand{\lie}[1]{\mathfrak{#1}}
\makeatletter
\def\section{\def\@secnumfont{\mdseries}\@startsection{section}{1}%
  \z@{.7\linespacing\@plus\linespacing}{.5\linespacing}%
  {\normalfont\scshape\centering}}
\def\subsection{\def\@secnumfont{\bfseries}\@startsection{subsection}{2}%
  {\parindent}{.5\linespacing\@plus.7\linespacing}{-.5em}%
  {\normalfont\bfseries}}
\makeatother
\def\subl#1{\subsection{}\label{#1}}
 \nc{\Hom}{\operatorname{Hom}}
\nc{\End}{\operatorname{End}} \nc{\wh}[1]{\widehat{#1}} \nc{\Ext}{\operatorname{Ext}} \nc{\ch}{\text{ch}} \nc{\ev}{\text{ev}}
\nc{\Ob}{\operatorname{Ob}} \nc{\soc}{\operatorname{soc}} \nc{\rad}{\operatorname{rad}} \nc{\head}{\operatorname{head}}
\def\Im{\operatorname{Im}}
\def\gr{\operatorname{gr}}
\def\mult{\operatorname{mult}}

\nc{\krsm}{KR^\sigma(m\omega_i)} \nc{\krsmzero}{KR^\sigma(m_0\omega_i)} \nc{\krsmone}{KR^\sigma(m_1\omega_i)}
 \nc{\vsim}{v^\sigma_{i,m}}

 \nc{\Cal}{\cal} \nc{\Xp}[1]{X^+(#1)} \nc{\Xm}[1]{X^-(#1)}
\nc{\on}{\operatorname} \nc{\Z}{{\bold Z}} \nc{\J}{{\cal J}} \nc{\C}{{\bold C}} \nc{\Q}{{\bold Q}}
\renewcommand{\P}{{\cal P}}
\nc{\N}{{\Bbb N}} \nc\boa{\bold a} \nc\bob{\bold b} \nc\boc{\bold c} \nc\bod{\bold d} \nc\boe{\bold e} \nc\bof{\bold f} \nc\bog{\bold g}
\nc\boh{\bold h} \nc\boi{\bold i} \nc\boj{\bold j} \nc\bok{\bold k} \nc\bol{\bold l} \nc\bom{\bold m} \nc\bon{\bold n} \nc\boo{\bold o}
\nc\bop{\bold p} \nc\boq{\bold q} \nc\bor{\bold r} \nc\bos{\bold s} \nc\bou{\bold u} \nc\bov{\bold v} \nc\bow{\bold w} \nc\boz{\bold z}
\nc\boy{\bold y} \nc\ba{\bold A} \nc\bb{\bold B} \nc\bc{\bold C} \nc\bd{\bold D} \nc\be{\bold E} \nc\bg{\bold G} \nc\bh{\bold H} \nc\bi{\bold I}
\nc\bj{\bold J} \nc\bk{\bold K} \nc\bl{\bold L} \nc\bm{\bold M} \nc\bn{\bold N} \nc\bo{\bold O} \nc\bp{\bold P} \nc\bq{\bold Q} \nc\br{\bold R}
\nc\bs{\bold S} \nc\bt{\bold T} \nc\bu{\bold U} \nc\bv{\bold V} \nc\bw{\bold W} \nc\bz{\bold Z} \nc\bx{\bold x} {\title[Highest weight
categories of representations]{ Current algebras, Highest weight categories and quivers}
\author{Vyjayanthi Chari and  Jacob Greenstein}
\thanks{This work was partially supported by the NSF grant DMS-0500751}
\address{Department of Mathematics, University of
California, Riverside, CA 92521.} \email{chari@math.ucr.edu}
 \email{jacob.greenstein@ucr.edu}\maketitle
\begin{abstract}
We study the category of graded finite-dimensional 
representations of the polynomial current algebra associated 
to a simple Lie algebra. We prove that the category has 
enough injectives and compute the graded character of the 
injective envelopes of the simple objects as well as 
extensions betweeen simple objects. The simple objects in 
the category are parametized by the affine weight lattice. 
We show that with respect to a suitable refinement of the 
standard ordering on affine the weight lattice the category 
is highest weight. We compute the $\Ext$ quiver of the algebra 
of endomorphisms of the injective cogenerator of the 
subcategory associated to an interval closed finite subset of 
the weight lattice. Finally, we prove that there is a large 
number of interesting quivers of finite, affine and tame 
type that arise from our study. We also prove that the path 
algebra of star shaped quivers are the $\Ext$-algebra of a 
suitable subcategory.
\end{abstract}

\section*{Introduction}
In this paper we study the category $\cal G$ of  graded
finite-dimensional representations  of the polynomial current
algebra~$\lie g[t]$ associated to a simple finite dimensional Lie algebra~$\lie
g$. There are numerous interesting and related families of examples
of such representations: the Demazure modules arising from the
positive level representations of the affine algebra, the fusion
products of finite-dimensional representations of $\lie g[t]$ defined
in \cite{FL}, the Kirillov-Reshetikhin modules studied in
\cite{CMkir1, CMkir2} and the Weyl modules introduced in~\cite{CP}
and studied in \cite{CL,FoL}. All these representations are in
general reducible but always indecomposable.

The isomorphism classes of simple objects in $\cal G$ are indexed by
the set $\Lambda= P^+\times \bz_+$ where $P^+$ is the set of
dominant integral weights of $\lie g$. The set $\Lambda$ can be
identified in a natural way with a subset of the lattice of integral
weights $\wh P$ of the untwisted affine Lie algebra associated to
$\lie g$. We define an interval finite partial order $\preccurlyeq$
on $\Lambda$ which is a refinement of the usual order on $\wh P$ and
show (\thmref{thmone}) that  $\cal G$ is a highest weight category,
in the sense of~\cite{CPS}, with the poset of weights
$(\Lambda,\preccurlyeq)$.  To do this, we study first the category
$\wh{\cal G}$ of graded $\lie g[t]$-modules with finite-dimensional
graded pieces. This category has enough projectives and the graded
character of the projective modules can be described explicitly.
Then, using a certain duality, we are able to show that the category
$\cal G$ has enough injectives  and
 we  compute the graded character of the injective
envelope of any simple object. We  then prove that $\cal G$ is a
directed highest weight category by computing the  extensions
between simple objects.

In Section~\ref{TMPALG} we study algebraic structures associated
with Serre subcategories of~$\cal G$. For  an interval closed subset
$\Gamma$ of $\Lambda$, let $\cal G[\Gamma]$ be the full subcategory
of $\cal G$ consisting of objects whose simple constituents are
parametrized by elements of $\Gamma$ and let  $I(\Gamma)_\Gamma$ be
the injective cogenerator of $\cal G[\Gamma]$.  It is well-known
that there is an equivalence of categories between $\cal G[\Gamma]$
and the category of finite-dimensional right $\mathfrak
A(\Gamma)=\End_{\cal G[\Gamma]}I(\Gamma)_\Gamma$-modules. Moreover $\mathfrak
A(\Gamma)$ is a quotient of the path algebra of its $\Ext$ quiver
$Q(\Gamma)$ and has a compatible grading. By  using the character
formula for the injective envelopes, we show that $Q(\Gamma)$  can
be computed quite explicitly in terms of finite dimensional
representations of~$\lie g$.

In Sections~\ref{EXH} and~\ref{QR} we show that there are many
interesting quivers arising from our study. Thus, in
Section~\ref{EXH} we see that for all $\lie g$ (in some cases one
has to exclude $\lie{sl}_2$ or~$\lie g$ of type~$C_\ell$), there
exists  interval closed finite subsets $\Gamma$ such that the
corresponding algebra $\mathfrak A(\Gamma)$ is hereditary and
$Q(\Gamma)$ is (a) a generalized Kronecker quiver; (b) a quiver of
type $\mathbb A_\ell$, $\mathbb D_\ell$; (c) an affine quiver of
type $\tilde{\mathbb D}_\ell$; (d) any star shaped quiver with three
branches. In Section~\ref{QR} we study an example which arises from
the theory of Kirillov-Reshetikhin modules for $\lie g[t]$ where
$\lie g$ is of type $D_n$. In this case the  algebra $\mathfrak
A(\Gamma)$ is not hereditary, but is still of tame representation
type.

\subsection*{Acknowledgements}
The first author is very grateful to Steffen Koenig for his infinite
patience in answering many questions and for his generosity in
providing references and detailed explanations, this paper could not
have been written without those discussions. The second author
thanks Olivier Schiffmann and Wolfgang Soergel. Part of this work
was done while the first author was visiting the University of
Cologne and the second author was visiting the Weizmann Institute of
Science. It is a pleasure to thank Peter Littelmann and the algebra
group of the University of Cologne and Anthony Joseph for their
hospitality. We also thank Brian Parshall for pointers to references
in the literature. Finally, we are grateful to Claus Michael Ringel
for explaining to us the proof that  the example in Section~\ref{QR}
is of tame type.

\section{The category~$\mathcal G$}\label{CAT}

\subsection{The simple Lie algebras and the associated current algebras}\label{CAT10}
Throughout the paper $\lie g$ denotes a finite-dimensional complex
simple Lie algebra and $\lie h$ a fixed Cartan subalgebra of $\lie
g$. Set~$I=\{1,\dots,\dim\lie h\}$ and let  $\{\alpha_i: i\in
I\}\subset\lie h^*$ be  a set of simple roots of $\lie g$ with
respect to $\lie h$. Let $R\subset\lie h^*$ (respectively, $R^+$,
$P^+$, $Q^+$) be the corresponding set of roots  (respectively,
positive roots, dominant integral weights, the $\bz_+$-span of
$R^+$) and let $\theta\in R^+$ be the highest root.
Let $W\subset \operatorname{Aut}(\lie h^*)$ be the Weyl group of $\lie g$ and
$w_\circ$ be the longest element of $W$. For
$\alpha\in R$ denote by $\lie g_\alpha$ the corresponding root
space. The subspaces $\lie n^\pm=\bigoplus_{\alpha\in R^+}\lie
g_{\pm\alpha},$ are Lie subalgebras of $\lie g$. Fix a Chevalley
basis $x^\pm_\alpha$, $\alpha\in R^+$, $h_i$, $i\in I$ of $\lie g$
and for~$\alpha\in R^+$, set~$h_\alpha=[x_\alpha,x_{-\alpha}]$. Note that $h_{\alpha_i}=h_i$, $i\in I$. For
$i\in I$, let $\omega_i\in P^+$ be defined by
$\omega_i(h_j)=\delta_{ij}$ for all $j\in I$.

Let $\cal F(\lie g)$ be the category of finite-dimensional $\lie g$-modules with the morphisms being maps of  $\lie g$-modules.
 In particular, we write $\Hom_{\lie g}$ for~$\Hom_{\cal
F(\lie g)}$. The set $P^+$ parametrizes the isomorphism classes of simple objects in~$\cal F(\lie g)$. For $\lambda\in P^+$, let $V(\lambda)$ be
the  simple module in  the corresponding isomorphism class which  is generated by an element $v_\lambda\in V(\lambda)$ satisfying the defining
relations:
$$
\lie n^+ v_\lambda=0,\quad hv_\lambda=\lambda(h)v_\lambda,\quad (x^-_{\alpha_i})^{\lambda(h_i)+1}v_\lambda =0,
$$
for all~$h\in\lie h$, $i\in I$. The module~$V(-w_\circ\lambda)$  is the $\lie g$-dual of $V(\lambda)$.
 If $V\in\cal F(\lie g)$, write
$$V=\bigoplus_{\lambda\in\lie h^*} V_\lambda,$$ where $V_\lambda=\{v\in
V:hv=\lambda(h)v,\ \ \forall\ h\in\lie h\}$. Set $\wt(V)=\{\lambda\in\lie h^*:V_\lambda\ne 0\}$. Finally, recall also that the category~$\cal
F(\lie g)$ is semi-simple, i.e. any object in~$\cal F(\lie g)$
is isomorphic to a direct sum of the modules $V(\lambda)$, $\lambda\in P^+$. We shall use
the following standard results in the course of the paper (cf.~\cite{PRV} for \eqref{CAT10.iv}).
\begin{lem}
\begin{enumerate}[{\rm(i)}]
\item\label{CAT10.i} Let $\lambda\in P^+$. Then
 $\wt(V(\lambda))\subset \lambda-Q^+$.
\item\label{CAT10.ii} Let $V\in\cal F(\lie g)$.Then $w\wt(V)\subset\wt(V)$ for all $w\in W$
 and $\dim V_\lambda=V_{w\lambda}$.
\item\label{CAT10.iii} Let $V\in\cal F(\lie g)$. Then $$\dim\Hom_{\lie
g}(V(\lambda), V)=\dim\{v\in V_\lambda: \lie n^+ v=0\}.$$
\item\label{CAT10.iv} Let $\lambda,\mu\in P^+$. Then the module $V(-w_\circ\lambda)\otimes V(\mu)$ is
generated as a $\bu(\lie g)$-module by the element $v=v_{-\lambda}\otimes v_\mu$ with defining relations:
$$
(x^+_{\alpha_i})^{\mu(h_i)+1}v=(x^-_{\alpha_i})^{\lambda(h_i)+1}v=0,\qquad
hv=(\mu-\lambda)(h)v,
$$
for all $i\in I$ and $h\in\lie h$.\qed
\end{enumerate}
\end{lem}

Given any Lie algebra $\lie a$ let $\lie a[t]=\lie a\tensor \bc[t]$ be the polynomial current algebra of~$\lie a$. Let~$\lie a[t]_+$ be the
Lie ideal~$\lie a\tensor t \bc[t]$. Both $\lie a[t]$ and $\lie a[t]_+$ are $\bz_+$-graded Lie algebras with the grading given by powers of $t$.
Let  $\bu(\lie a)$ denote the universal enveloping algebra of $\lie a$. Then  $\bu(\lie a[t])$ has a natural $\bz_+$-grading as an associative
algebra and we let $\bu(\lie a)[k]$ be the $k^{th}$-graded piece.  The algebra $U(\lie a)$ is a Hopf algebra, the comultiplication being given
by extending the assignment $x\to x\otimes 1+1\otimes x$ for $x\in\lie a$ to an algebra homomorphism of $U(\lie a)$. In the case of $\bu(\lie
a[t])$ the  comultiplication is a map of graded algebras.

{\em  In the course of the paper, we shall repeatedly use the fact that $\bu(\lie a[t])$ is generated as a graded algebra by $\lie a$ and $\lie
a\otimes t$ without further comment.}

\subsection{The category $\widehat{\cal G}$ }\label{CAT30}
Let $\wh{\cal G}$ be the category whose objects are graded $\lie g[t]$-modules with finite-dimensional graded pieces and where the morphisms are
graded maps of $\lie g[t]$-modules. More precisely, if $V\in\Ob\wh{\cal G}$ then
$$
V=\bigoplus_{r\in\bz_+} V[r],
$$
where $V[r]$ is a finite-dimensional subspace of $V$ such that $ (xt^k)V[r]\subset V[r+k] $ for all~$x\in\lie g$ and~$r,k\in\bz_+$. In
particular, $V[r]\in\Ob\cal F(\lie g)$ Also, if $V,W\in\Ob\wh{\cal G}$,then
$$\Hom_{\wh{\cal G}}(V,W)=\{f\in\Hom_{\lie g[t]}(V,W):
f(V[r])\subset W([r],\, r\in\bz_+\}.$$ For $f\in\Hom_{\wh{\cal G}}(V,W)$ let $f[r]$ be the  restriction of $f$ to $V[r]$. Clearly,
$f[r]\in\Hom_{\lie g}(V[r],W[r])$.

Define a covariant functor  $\ev:\cal F(\lie g)\to \wh{\cal G}$ by the requirements:
$$\ev(V)[0]=V,\qquad \ev(V)[r]=0,\qquad r>0,$$
and with $\lie g[t]$-action given by
$$
(xt^k)v=\delta_{k,0} xv,\qquad \forall\,x\in\lie g,\, k\in\bz_+,\, v\in V
$$
and
\begin{equation}\label{CAT30.10}
\Hom_{\wh{\cal G}}(\ev(V),\ev(W))= \Hom_{\lie g}(V,W).
\end{equation}
For~$s\in\bz_+$ let~$\tau_s$ be the grading shift given by
$$
(\tau_s V)[k]=V[k-s],\qquad k\in\bz_+, \quad V\in\Ob{\wh{\cal G}}.
$$
Clearly $\tau_s(V)\in\Ob{\wh{\cal G}}$.

\subsection{Simple objects in $\wh{\cal G}$}\label{CAT40}
For $\lambda\in P^+$ and $r\in\bz_+$, set
\begin{equation}\label{CAT40.20}
V(\lambda,r)=\tau_r(\ev(V(\lambda)).
\end{equation}
\begin{prop}The isomorphism classes of simple objects in $\wh{\cal G}$
are parametrized by pairs $(\lambda,r)$ and we have
\begin{align*}
&\Hom_{\wh{\cal G}}(V(\lambda,r), V(\mu,s))=0,\qquad (\lambda,r)\ne (\mu,s),\\
&\Hom_{\wh{\cal G}}(V(\lambda,r),V(\lambda,r))\cong\bc.
\end{align*}
 Moreover if
$V\in\Ob\cal G$ is such that $V=V[n]$ for some $n\in\bz_+$, then $V$ is semi-simple.
\end{prop}
\begin{pf}The modules $V(\lambda,r)$, $(\lambda,r)\in\Lambda$  are obviously simple and
non-isomorphic. Moreover, if $V\in\Ob\wh{\cal G}$ is such that $V\ne V[r]$ for some $r$, there exists $m, m'\in\bz_+$ with
$m'>m$ such that ~$V[j]\not=0$, for $j\in\{m,m'\}$.
Hence the subspace $\bigoplus_{k>m}V[k]$ is a nontrivial proper graded $\lie g[t]$-submodule of
$V$ and it follows that $V$ is not simple. Assume now that  $V$ is simple so that  $V=V[r]$ for some $r$. This implies that $V$ is
finite-dimensional and also that $\lie g[t]_+ V=0$.  It follows that $V$ must be isomorphic to $V(\lambda)$ for some $\lambda\in P^+$ as a $\lie
g$-module and hence $V\cong V(\lambda,r)$ as $\lie g[t]$-modules. The other statements  are now obvious.
\end{pf}

\subsection{Tensor structure of the category~$\wh{\cal  G}$}\label{CAT50}
Let~$V,W\in\Ob\wh{\cal G}$. Then~$V\tensor W$ is a $\lie g[t]$-module with the action being given by the comultiplication.
Given~$k\in\bz_+$, set
$$
(V\tensor W)[k]=\bigoplus_{i\in\bz_+} V[i]\tensor W[k-i],
$$
with the usual convention that~$W[j]=0$ if~$j<0$. The following is trivially checked.
\begin{lem}
\begin{enumerate}[{\rm(i)}]
\item\label{CAT50.i} $V\tensor W=\bigoplus_{k\in\bz_+} (V\tensor
W)[k]$ and for all $r\in\bz_+$, we have
$$
(xt^r)((V\tensor W)[k])\subset (V\tensor W)[k+r].
$$
In particular, $\wh{\cal G}$ is a tensor category.
\item\label{CAT50.ii} For all~$r,s\in\bz_+$
\begin{equation} \label{CAT50.10}
\tau_s V\cong V\tensor V(0,s),\qquad \tau_{r+s}(V\tensor W)\cong(\tau_r V)\tensor (\tau_s W).
\end{equation}
\end{enumerate}\qedhere
\end{lem}

\subsection{The subcategories $\cal G$ and $\cal G_{\le s}$} \label{CAT70}
Let $\cal G_{\le s}$ be the full subcategory of $\wh{\cal G}$ whose objects~$V$ satisfy
$$V[r]=0,\qquad \forall\, r>s,
$$
and let $\cal G$ be the full subcategory of $\cal G$ consisting of $V\in\Ob\wh{\cal G}$ such that $V\in\Ob\cal G_{\le s}$ for some $s\in\bz_+$.
It follows from the definition that $\cal G_{\le s}$ is a full subcategory of~$\cal G_{\le r}$ for all~$s< r\in\bz_+$.
Given $V\in\Ob\cal G$, let $\soc(V)\in\Ob\cal G$ be the maximal semi-simple subobject
of~$V$. Similalry, given $V\in\Ob\wh{\cal G}$, let~$\head(V)$ be the maximal semi-simple quotient of~$V$.

Given $s\in\bz_+$ and $V\in\Ob\wh{\cal G}$, define
$$
V_{>s}= \bigoplus_{r>s} V[r],\qquad V_{\le s}= V/V_{>s}
$$
Then~$V_{\le s}\in\Ob\cal G_{\le s}$.
 Furthermore, if $f\in\Hom_{\wh{\cal G}}(V,W)$, then $V_{>s}$ is contained in the kernel of  the canonical
morphism $\bar f:V\to W_{\le s}$ and hence we have a natural morphism~$f_{\le s}\in\Hom_{\wh{\cal G}_{\le s}}( V_{\le s}, W_{\le s})$.

\begin{lem}
\begin{enumerate}[{\rm(i)}]
\item\label{CAT70.i} For all $r,s\in\bz_+$, and $V\in\Ob\cal
G_{\le r}$, $W\in\Ob\cal G_{\le s}$ we have $$V\otimes W\in\Ob\cal G_{\le r+s}.$$ In particular $\cal G$ is a tensor subcategory of $\wh{\cal
G}$. \item\label{CAT70.ii} The assignments $V\mapsto V_{\le r}$ for all~$V\in\Ob\wh{\cal G}$ and~$f\mapsto f_{\le r}$ for
all~$f\in\Hom_{\wh{\cal G}}(V,W)$, $V,W\in\Ob\wh{\cal G}$ define a full, exact and essentially surjective functor from $\wh{\cal G}$ to $\cal G_{\le
r}$. \item\label{CAT70.iii} For any $V\in\Ob\wh{\cal G}$, $\lambda\in P^+$ and $r,s\in\bz_+$ with $s\ge r$, we have
$$(V\otimes V(\lambda,r))_{\le s}\cong V_{\le s-r}\otimes V(\lambda,r).$$
\end{enumerate}
\end{lem}
\begin{pf} Parts~\eqref{CAT70.i} and~\eqref{CAT70.ii} are obvious. For the last part, consider the natural map of graded $\lie g[t]$-modules
$V\otimes V(\lambda,r)\to V_{\le s-r}\otimes V(\lambda,r)$. The assertion follows by noting that  $(V\otimes V(\lambda,r))[k]=V[k-r]\otimes
V(\lambda,r)$ for all $k\in\bz_+$.
\end{pf}
From now on, given~$V\in\Ob\cal G$ we denote by~$[V:V(\lambda,r)]$ the multiplicity of~$V(\lambda,r)$ in a composition series
for~$V$. Furthermore, given~$W\in\Ob\wh{\cal G}$, we set~$[W:V(\lambda,r)]:=[W_{\le r}:V(\lambda,r)]$. Observe that
$[V:V(\lambda,r)]$ equals the $\lie g$-module multiplicity of~$V(\lambda)$ in~$V[r]$.
For any~$V\in\Ob\wh{\cal G}$, define
$$\Lambda(V)=\{(\lambda,r)\in\Lambda: [V:V(\lambda,r)]\ne 0\}.$$

\subsection{}\label{HWCAT}
We recall the following definition (which motivated much of this paper) of a directed  category following~\cite{CPS,PSW}, in the context of
interest to us. Thus let~$\cal C$ be an abelian category over~$\bc$ whose objects are complex vector spaces, have finite length and such that
$\Hom_{\cal C}(M,N)$ is finite-dimensional for all $M,N\in\Ob\cal C$.
\begin{defn} We say
that $\cal C$ is a directed category if
\begin{enumerate}[$1^\circ.$]
\item\label{HWCAT90.1} The simple objects in $\cal C$ are
parametrized by a poset $(\Pi,\le)$ ({\em the poset of weights}) such that for all~$\tau\in\Pi$, the set~$\{\xi\in\Pi\,:\,\xi<\tau\}$ is finite.
\item\label{HWCAT90.2} Given~$\sigma\in\Pi$, let~$S(\sigma)$ be a simple object in the corresponding isomorphism class. Then
$$
\Ext^1_{\cal C}(S(\sigma),S(\tau))\not=0 \implies \sigma<\tau.
$$
\end{enumerate}\end{defn}
It is immediate from~\cite{DR,Rin,Soe} that a directed category has enough injectives.
 Given~$\Xi\subset\Pi$,
let~$\mathcal C[\Xi]$ be the full subcategory of~$\cal C$ whose objects satisfy
$$
M\in\Ob\cal C[\Xi],\quad [M:S(\tau)]\not=0\implies \tau\in\Xi.
$$ It is clear that if $\cal C$ is a directed category with poset of weights $\Pi$, then  for any subset $\Xi\subset\Pi$,
the category  $\cal C[\Xi]$ is also directed with poset of weights $\Xi$.
Given $\sigma,\tau\in\Xi $ with $\sigma<\tau$,
we denote by $[\sigma,\tau]$ the interval $\{\xi:\sigma\le \xi\le \tau\}$. A subset $\Xi$ of $\Pi$ is said to be {\em interval closed} if
$\sigma<\tau\in\Xi$ implies that $[\sigma,\tau]\subset \Xi$.

\subsection{}\label{CAT110} We can now state the main result of this section. Set $$\Lambda=P^+\times\bz^+,\qquad \Lambda_{\le
r}=\{(\nu,l)\in\Lambda\,:\, l\le r\},\qquad r\in\bz_+.
$$
Define a strict partial order on~$\Lambda$ in the following way.
Given~$(\lambda,r), (\mu,s)\in\Lambda$, say that $(\mu,s)$ covers
$(\lambda,r)$ if and only if $s=r+1$ and~$\mu-\lambda\in
R\sqcup\{0\}$. It follows immediately that for
any~$(\mu,s)\in\Lambda$ the set of~$(\lambda,r)\in\Lambda$ such
that~$(\mu,s)$ covers $(\lambda,r)$ is finite. Let $\preccurlyeq$ be
the unique partial order on $\Lambda$ generated by this cover
relation. Then~$\{(\mu,s): (\mu,s)\prec(\lambda,r)\}$ is finite for
all~$(\lambda,r)\in\Lambda$. Note that if $(\lambda,r)\prec(\mu,s)\in \Lambda_{\le k}$
then $(-w_\circ\mu,k-s)\prec (-w_\circ\lambda,k-r)$.

 \begin{thm} \label{thmone} For all $\Gamma\subset\Lambda$,  the category $\cal G[\Gamma]$ is a directed category with poset
 of weights $(\Gamma,\preccurlyeq)$.
 \end{thm}
\noindent
We prove this theorem in the next section (\propref{CAT170}).

\section{Injective and projective objects}

 \subsection{Projectives in $\wh{\cal G}$ and $\cal G_{\le r}$}\label{CAT130}
Given $(\lambda,r)\in \Lambda$, set
$$
P(\lambda,r)=\bu(\lie g[t])\otimes_{\bu(\lie g)} V(\lambda,r) .
$$
Clearly, $P(\lambda,r)$ is an infinite dimensional graded $\lie g[t]$-module. Using  the PBW theorem  we have an isomorphism of graded vector
spaces
$$
\bu(\lie g[t])\cong \bu(\lie g[t]_+)\tensor \bu(\lie g),
$$
and hence we get
\begin{equation}\label{CAT130.10}
P(\lambda,r)[k]=\bu(\lie g[t]_+)[k-r]\otimes V(\lambda,r),
\end{equation}
where we understand that $\bu(\lie g[t]_+)[k-r]=0$ if $k<r$.  This shows that  $P(\lambda,r)\in\Ob\wh{\cal G}$ and also that
$$P(\lambda,r)[r] =1\otimes V(\lambda,r).$$ Set
$p_{\lambda,r}=1\otimes v_{\lambda,r}$.
\begin{prop}
Let $(\lambda,r)\in\Lambda$,  $s\in\bz_+$ and $s\ge r$.
\begin{enumerate}[{\rm(i)}]
\item\label{CAT130.i} $P(\lambda,r)$ is  generated as a $\lie
g[t]$-module by $p_{\lambda,r}$ with defining relations:
$$(\lie n^+)p_{\lambda,r}=0,\quad
hp_{\lambda,r}=\lambda(h)p_{\lambda,r},\quad (x^-_{\alpha_i})^{\lambda(h_i)+1}p_{\lambda,r}=0,$$ for all $h\in\lie h$, $i\in I$. Hence,
$P(\lambda,r)$ is the projective cover in the category~$\wh{\cal G}$ of its unique simple quotient $V(\lambda,r)$. Moreover the kernel
$K(\lambda,r )$ of the canonical projection $P(\lambda,r)\twoheadrightarrow V(\lambda,r)$ is generated as a $\lie g[t]$-module by $P(\lambda,r)[r+1]=\lie
g\otimes V(\lambda,r)$.

\item\label{CAT130.ii}
$P(\lambda,r)\cong P(0,0)\otimes V(\lambda,r)$ as objects in~$\wh{\cal G}$.
\item\label{CAT130.iii} The modules $P(\lambda,r)_{\le s}$ are
projective in $\cal G_{\le s}$ and
\begin{equation*}
P(\lambda,r)_{\le s}\cong P(0,0)_{\le s-r}\otimes V(\lambda,r)
\end{equation*}
\item\label{CAT130.iv} As $\lie g$-modules, we have
$$P(0,0)[k]\cong\bu(\lie g[t]_+)[k]\cong \bigoplus_{(r_1,\dots,r_k)\in\bz_+^k\,:\,
 \sum_{j=1}^kjr_j=k}S^{r_1}(\lie g)\tensor\cdots\tensor S^{r_k}(\lie g),
$$ where
 $S^p(\lie g)$ denotes the $p^{th}$ symmetric power
of $\lie g$.
\item\label{CAT130.v}  Let~$(\mu,s)\in\Lambda$. Then~$[K(\lambda,r):V(\mu,s)]\not=0$ only if~$(\lambda,r)\prec(\mu,s)$.
\item\label{CAT130.vi}  Let~$V\in\Ob\wh{\cal G}$. Then~$\dim\Hom_{\wh{\cal G}}(P(\lambda,r),V)=[V:V(\lambda,r)]$.
\end{enumerate}
\end{prop}
\begin{pf} The fact that $P(\lambda,r)$ is projective  in the category~$\wh{\cal G}$ is standard in
relative homological algebra (cf.~\cite{Hoch}). The other statements in~\eqref{CAT130.i} are immediate from the discussion preceding the
proposition. For part~\eqref{CAT130.ii}, note that the element $p_{0,0}\otimes v_{\lambda,r}$ satisfies the defining relations of
$P(\lambda,r)$. Moreover it is easily seen that
$$
\bu(\lie g[t])(p_{0,0}\otimes v_{\lambda,r})=P(0,0)\otimes V(\lambda,r).
$$
Hence we have a surjective morphism $P(\lambda,r)\to P(0,0)\otimes V(\lambda,r)$ in~$\wh{\cal G}$. On the other hand, \eqref{CAT130.10} implies
that $P(\lambda,r)\cong P(0,0)\otimes V(\lambda,r)$ as vector spaces and~\eqref{CAT130.ii} is proved. Part~\eqref{CAT130.iii} is immediate from
\lemref{CAT70}(\ref{CAT70.ii},\ref{CAT70.iii}). The first isomorphism in \eqref{CAT130.iv} is obvious. To prove the second, we may
assume that $k>0$ since $\bu(\lie g[t]_+)[0]=\bc$. For $s\ge 0$, let $\bu(\lie g[t]_+)_{\le s}$, $s\ge 0$ be the subspace
of~$\bu(\lie g[t]_+)$ spanned by the set
$$\{(y_1t^{r_1})\cdots (y_kt^{r_k}): y_jt^{r_j}\in\lie g[t]_+,\,1\le j \le k\le s\},$$ and set
$$
\bu(\lie g[t]_+)[k]_{\le s}=\bu(\lie g[t]_+)[k]\cap \bu(\lie g[t]_+)_{\le s}.
$$
Then $\bu(\lie g[t]_+)[k]=\bu(\lie g[t]_+)[k]_{\le k}$ and
for $0\le s\le k$ the subspaces  $\bu(\lie g[t]_+)[k]_{\le s}$ define an increasing filtration on $\bu(\lie g[t]_+)[k]$. Moreover, regarding
$\lie g[t]_+$ as a $\lie g$-module via the adjoint action on~$\lie g[t]$, we see that this filtration is in fact $\lie g$-equivariant.
 Since $\bu(\lie g[t]_+)[k]_{\le r}$
is finite dimensional we get an isomorphism of $\lie g$-modules,
\begin{align*}
\bu(\lie g[t]_+)[k]_{\le r}&\cong_{\lie g} \bu(\lie g[t]_+)[k]_{\le r-1}\oplus (\bu(\lie g[t]_+)[k]_{\le r}/\bu(\lie g[t]_+)[k]_{\le r-1})
\\&\cong_{\lie g} \bu(\lie g[t]_+)[k]_{\le r-1} \oplus S^r(\lie g[t]_+)[k],
\end{align*}
the second isomorphism being a consequence of the PBW theorem.
It follows by a downward induction on~$r$ that
$$\bu(\lie g[t]_+)[k]\cong_{\lie g}
\bigoplus_{r=1}^k S^r(\lie g[t])_+[k]=S(\lie g[t]_+)[k].
$$
Given  a partition $\bon= (n_s\ge n_{s-1}\ge \cdots \ge n_1>0)$ of $k$, let  $V(\bon)$ be the subspace of~$S(\lie g[t]_+)[k]$ spanned by the
elements $$ \{(x_{j_1} t^{n_1})\cdots (x_{j_s} t^{n_s}): x_j\in\lie g\}.
 $$
Clearly $V(\bon)$ is  a $\lie g$--submodule of $S(\lie g[t]_+)[k]$ and we have  $$
  S(\lie g[t]_+)[k]=\bigoplus_{\bon\vdash k}
V(\bon).
$$
The result follows since $$V(\bon)\cong_{\lie g} S^{r_1}(\lie g)\tensor \cdots\tensor S^{r_k}(\lie g),$$ where $r_j=|\{ 1\le r\le s\,:\,
n_r=j\}|$.

Part~\eqref{CAT130.v} is now obvious. To establish \eqref{CAT130.vi}, it is enough to observe that
by~\eqref{CAT130.i} the natural map
$$\Hom_{\wh{\cal G}}(P(\lambda,r),V)\to\Hom_{\cal G_{\le r}}(P(\lambda,r)_{\le r},V_{\le r})$$
is injective and hence is an isomorphism (\propref{CAT70}\eqref{CAT70.ii}). The statement follows since $P(\lambda,r)_{\le r}\cong V(\lambda,r)$ is projective
in~$\cal G_{\le r}$ and every object in~$\cal G_{\le r}$ has finite length.
\end{pf}
In what follows, we shall write
$$S^{(k)}(\lie g)=\bigoplus_{(r_1,\dots,r_k)\in\bz_+^k\,:\,
 \sum_{j=1}^kjr_j=k} S^{r_1}(\lie g)\tensor\cdots\tensor S^{r_k}(\lie
 g).
$$
Observe that~$S^{(k)}(\lie g)$ is a $\lie g$-module quotient of~$\lie g^{\tensor k}$. Indeed, the map~$\lie g^{\tensor k}\to
\bu(\lie g[t]_+)[k]\cong_{\lie g} S^{(k)}(\lie g)$ given by extending~$x_1\tensor\cdots\tensor x_k\mapsto (x_1 t)\cdots(x_k t)$, $x_j\in\lie g$, $1\le j\le k$ is a
surjective $\lie g$-module homomorphism.

\subsection{Morphisms between projectives}\label{CAT150} The next
proposition is trivial, but we include it here explicitly since it is used repeatedly in the later sections where we construct examples of
interesting subcategories of $\cal G$.

Let~$s\le r\le \ell\in\bz_+$. It is immediate from the definitions and Proposition~\ref{CAT130} that
\begin{equation}\label{CAT150.10}
\begin{split}
\Hom_{\cal G}(P(\lambda,r)_{\le \ell},P(\mu,s)_{\le
\ell})&=\Hom_{\wh{\cal G}}(P(\lambda,r),P(\mu,s))\\
&\cong \Hom_{\lie g}(V(\lambda),\bu(\lie g[t]_+)[r-s]\tensor V(\mu)).
\end{split}
\end{equation}
In this section we make the last isomorphism explicit. Let $s\le r\in\bz_+$ and $\lambda,\mu\in P^+$.  Given $f\in\Hom_{\lie
g}(V(\lambda),\bu(\lie g[t]_+)[r-s]\tensor V(\mu))$, define $\bof\in\Hom_{\bc}(P(\lambda,r),P(\mu,s))$ by
$$
\bof(u\tensor v)=\sum_p uu_p\tensor v_p,
$$
where~$u\in \bu(\lie g[t]_+)$, $v\in V(\lambda,r)$ and~$\sum_p u_p \tensor v_p=f(v)$. Let $m:\bu(\lie g[t])\otimes \bu(\lie g[t])\to\bu(\lie
g[t])$ be the multiplication map. The next proposition is a straightforward consequence of Proposition~\ref{CAT130} and we omit the details of
the calculations.
\begin{prop} Let $(\lambda,r), (\mu,s)\in\Lambda$.
\begin{enumerate}[{\rm(i)}]
\item\label{CAT150.i} The assignment~$\bof\mapsto \bof[r]$ defines
an isomorphism of vector spaces
$$
\phi^{(\lambda,r)}_{(\mu,s)}: \Hom_{\wh{\cal G}} (P(\lambda,r),P(\mu,s))\to \Hom_{\lie g}(V(\lambda),\bu(\lie g[t]_+)[r-s]\tensor V(\mu)).
$$
Moreover if $\bof\in\Hom_{\wh{\cal G}}(P(\lambda,r), P(\mu,s))$, $\bog\in\Hom_{\wh{\cal G}}(P(\mu,s), P(\nu,k))$, then
$$(\bog\circ\bof)[r]=(m\otimes 1)\circ(1\otimes \bog[s])\circ \bof[r].$$

\item\label{CAT150.ii} The assignment~$f\mapsto \bof$ is an
isomorphism
$$
\psi^{(\lambda,r)}_{(\mu,s)}: \Hom_{\lie g}(V(\lambda),\bu(\lie g[t]_+)[r-s]\tensor V(\mu))\to \Hom_{\wh{\cal G}}(P(\lambda,r),P(\mu,s))
$$
which is the inverse of~$\phi^{(\lambda,r)}_{(\mu,s)}$. Moreover, if~$f\in \Hom_{\lie g}(V(\lambda),\bu(\lie g[t]_+)[r-s]\tensor V(\mu))$,
$g\in\Hom_{\lie g}(V(\mu),\bu(\lie g[t]_+)[s-k]\tensor V(\nu))$ then
$$
\psi^{(\lambda,r)}_{(\nu,k)}((m\otimes 1)\circ(1\tensor g)\circ f)= \psi^{(\mu,s)}_{(\nu,k)}(g)\circ \psi^{(\lambda,r)}_{(\mu,s)}(f).
$$
\end{enumerate}
\end{prop}

\subsection{Duality in~$\cal G_{\le s}$}\label{CAT90}
Given~$V,W\in\Ob\cal G$, the vector space  $\Hom_{\bc}(V,W)$ is a $\lie g[t]$-module with respect to the usual action
$$
((xt^r)\cdot f)(v)= (xt^r) f(v)-f( (xt^r) v),
$$
for all~$x\in\lie g$, $r\in\bz_+$, $f\in\Hom_{\bc}(V,W)$ and~$v\in V$. For $k\in\bz_+$, set
\begin{equation}\label{CAT90.10}
\Hom_{\bc}(V,W)[k]=\bigoplus_{i\in\bz_+} \Hom_{\bc}(V[i],W[i+k])
\end{equation}
and define
$$
\Hom_{\bc}^+(V,W)=\bigoplus_{k\in\bz_+} \Hom_{\bc}(V,W)[k].
$$
Since~$V[i]=0$ for all but finitely many~$i\in\bz_+$,
$\dim\Hom_{\bc}(V,W)[k]<\infty$. Notice that
$\Hom_{\bc}^+(V,W)=\Hom_{\bc}(V,W)$ provided that~$V\in\Ob\cal
G_{\le r}$ and~$W[i]=0$ for all~$i<r$. The proof of the following
proposition is quite standard and is omitted.
\begin{prop} Let $V,W\in\Ob\cal G$, $r,s\in\bz_+$.
\begin{enumerate}[{\rm(i)}]
\item\label{CAT90.0} For all~$x\in\lie g$, $i,r,k\in\bz_+$ and
$f\in\Hom_{\bc}(V[i],W[i+k])$,
\begin{equation*}
(xt^r)\cdot f\in\Hom_{\bc}(V[i],W[i+k+r])\oplus \Hom_{\bc}(V[i-r],W[i+k]).
\end{equation*}
In particular, $\Hom_\bc^+(V,W)\in\Ob\cal G$ and
$$
W\in\Ob\cal G_{\le s}\implies \Hom_{\bc}^+(V,W)\in\Ob\cal G_{\le s}.
$$
\item\label{CAT90.i} Let ${}^{\#_s}:\cal G_{\le s}\to \cal G_{\le
s}$ be the contravariant functor defined by $\Hom_{\bc}(-,V(0,s))$. Then
 ${}^{\#_s}$ is exact and for all $\lambda\in P^+$, $r\le s$,
$$
(V^{\#_r})^{\#_s}\cong \tau_{s-r} V,\qquad V(\lambda,r)^{\#_s}\cong V(-w_\circ\lambda,s-r).$$ In particular, ${}^{\#_s}$ defines an involutive
auto-duality on the category~$\cal G_{\le s}$.

\item\label{CAT90.ii} Suppose that~$V\in\Ob \cal G_{\le r}$,
$W\in\Ob\cal G_{\le s}$. As objects in $\cal G$,
$$
(V\otimes W)^{\#_{r+s}}\cong V^{\#_r}\otimes W^{\#_s},\qquad V\tensor W^{\#_s}\cong\Hom_{\bc}(W,V(0,s)\tensor V).
$$
\item\label{CAT90.iii} For all $V,W'\in\Ob\cal G$, $W\in\Ob\cal G_{\le s}$, we have an
isomorphism of vector spaces,
\begin{equation*}\tag*{\qedsymbol}\Hom_{\cal G}(V,W\tensor W')\cong
\Hom_{\cal G}(V\tensor W^{\#_s},V(0,s)\tensor W').
\end{equation*}
\end{enumerate}
\end{prop}

\subsection{Injective objects in $\cal G$ and ${\wh{\cal G}}$}\label{injG}
We begin with the following remark: any injective object of $\cal G$
is also injective in $\wh{\cal G}$. To prove this, let $I\in\Ob\cal G$
be injective and assume that $I[s]=0$ for all $s\ge r$. Suppose that
$\iota\in\Hom_{\wh{\cal G}}(V,W)$ is injective and
let $f\in\Hom_{\wh{\cal G}}(V,I)$. Since $f_{\le r}\in\Hom_{\cal
G}(V_{\le r}, I)$ and $\iota_{\le r}\in\Hom_{\wh{\cal G}}(V_{\le
r},W_{\le r})$ is injective there exists $\tilde f\in\Hom_{\cal G}(
W_{\le r},I)$ such that $\tilde f\circ\iota_{\le r}=f_{\le r}$. Let
$\bar f=\tilde f\circ p_r(W)$ where $p_r(W):W\to W_{\le r}$ be the
canonical projection. It is now easily checked that $\bar
f\circ\iota=f$.

For~$(\lambda,r)\in\Lambda$, set
$$ I(\lambda,r)\cong P(-w_\circ\lambda,0)_{\le r}{}^{\#_r}.
$$ 
It follows from~\propref{CAT170} and~\propref{CAT90} that $I(\lambda,r)\cong I(\lambda,0)\tensor V(0,r)$.
\begin{prop}
Let~$(\lambda,r)\in\Lambda$.
\begin{enumerate}[{\rm(i)}]
\item\label{injG.i}  The object $I(\lambda,r)$ is the injective envelope of $V(\lambda,r)$ in $\cal
G$.
\item\label{injG.ii}
For $k\in\bz_+$ we have
$$
I(\lambda,r)[r-k]\cong_{\lie g} S^{(k)}(\lie g)\tensor V(\lambda).
$$
\item\label{injG.iii}
Let~$(\mu,s)\in\Lambda$. Then~$[I(\lambda,r)/V(\lambda,r):V(\mu,s)]\not=0$ only if~$(\mu,s)\prec(\lambda,r)$.
\end{enumerate}
\end{prop}
\begin{pf} It is immediate from \propref{CAT90} that
$I(\lambda,r)$ is injective in $\cal G_{\le r}$. To prove that
$I(\lambda,r)$ is injective in $\cal G$, it suffices now to show
that $\Ext_{\cal G}^1(V(\mu,s),I(\lambda,r))=0$ if $s>r$, in other words that
every short exact sequence of the form $$0\to I(\lambda,r)\to V\to
V(\mu, s)\to 0$$ splits if $s>r$. Writing
$V=\bigoplus_{k\in\bz_+}V[k]$, we see that $V[k]=0$ if $k>s$. Hence
$\lie g[t]V[s]\subset V[s]$. Moreover, since $\bigoplus_{k\le r}
V[k]\cong_{\cal G} I(\lambda,r)$ it follows now that we have a
decomposition of $\lie g[t]$-modules $$V\cong I(\lambda,r)\oplus
V(\mu,s)$$ and hence the short exact sequence splits. To prove that
it is the injective envelope of $V(\lambda,r)$ it suffices to use
\propref{CAT130} and \lemref{CAT10} to notice that $V(\lambda,r)$ is
the unique irreducible subobject of $I(\lambda,r)$.

 The proof of
\eqref{injG.ii} and \eqref{injG.iii} are immediate from
\lemref{CAT10} and \propref{CAT130}.
\end{pf}
\begin{cor}\label{injG.cor}  Let $V,W\in\Ob\cal G$. 
\begin{enumerate}[{\rm(i)}]
\item\label{injG.cor.i} For all $j\ge 0$, we have
$$\Ext^j_{\cal G}(V,W)\cong\Ext^j_{\wh{\cal G}}(V,W).$$
\item\label{injG.cor.ii} Let~$I$ be the injective envelope of~$V$. If $(\lambda,r)\in \Lambda(I)$ then
$(\lambda,r)\preccurlyeq (\mu,s)$ for some~$(\mu,s)\in\Lambda(\soc(V))$.
\end{enumerate}
\end{cor}

\subsection{Extensions between simple objects} \label{CAT170}
The following proposition  proves~\thmref{thmone}.
\begin{prop}
For $(\lambda,r),(\mu,s)\in\Lambda$, we have
\begin{align*} &\Ext^1_{{\cal G}}(V(\lambda,r),V(\mu,s))=0,
\qquad s\ne r+1,\\
&\Ext^1_{{\cal G}}(V(\lambda,r),V(\mu,r+1)))\cong\Hom_{\lie g}(
V(\lambda), \lie g\otimes V(\mu)).
\end{align*}
In other words, $\Ext^1_{{\cal G}}(V(\lambda,r),V(\mu,s))=0$
unless $(\mu,s)$ covers $(\lambda,r)$.
\end{prop}
\begin{pf}
Applying $\Hom_{\cal G}(V(\lambda,r),-)$ to
the short exact sequence $$0\to V(\mu,s)\to
I(\mu,s)\to J(\mu,s)\to 0$$
gives
$$\Hom_{{\cal G}}(V(\lambda, r),
J(\mu,s))\cong\Ext^1_{\cal G}(V(\lambda,r), V(\mu,s)).
$$
The proposition obviously follows if we prove that
$$
\Hom_{{\cal G}}(V(\lambda, r),J(\mu,s))\cong\begin{cases} \Hom_{\lie g}(\lie g\otimes V(\lambda),V(\mu)),&\text{if $s=r+1$},\\ 0,&
\text{otherwise}.\end{cases}
$$
Let $\psi: V(\lambda,r)\to
J(\mu,s)$ be a non-zero element of $\Hom_{{\cal G}}(V(\lambda, r),
J(\mu,s))$. It follows from Proposition \ref{injG} that
$(\lambda,r)\prec (\mu,s)$ and hence in particular that $r<s$.
Suppose that $r<s-1$. Since~$V(\mu,s)$ is essential in~$I(\mu,s)$, there exists $V\subset I(\mu,s)$ such that
$V/V(\mu,s)\cong \psi(V(\lambda,r))$. Then $V=V[s]\oplus V[r]$ and
since $r<s-1$, it follows that $\lie g[t]V[r]\subset V$. Hence
$V[r]$ is in $\soc(I(\mu,s))$ which is impossible. Thus, $s=r-1$.

The following isomorphisms  which are consequences of \propref{injG}
establish the proposition.
\begin{gather*}\Hom_{\lie g}( V(\lambda),
\lie g\otimes V(\mu))\cong \Hom_{\cal G}(V(\lambda,r),
\tau_{r}\ev(\lie g\otimes V(\mu))),\\
\tau_r\ev(J(\mu,r+1)[r])\cong_{\cal G}
\tau_{r}\ev(\lie g\otimes V(\mu)),\\
\Hom_{{\cal G}}(V(\lambda,r),J(\mu,r+1)) \cong \Hom_{\cal
G}(V(\lambda,r),\tau_r\ev(J(\mu,r+1)[r])).\qedhere
\end{gather*}
\end{pf}

\subsection{}\label{CAT175}

Given $\Gamma\subset\Lambda$, set
\begin{gather*}{V_\Gamma}^+= \{v\in
V[s]_\mu: \lie n^+ v=0,\,\, (\mu,s)\in \Gamma\},\\
V_\Gamma=\bu(\lie g)V_\Gamma{}^+,\quad V^\Gamma=V/V_{\Lambda\setminus\Gamma}.
\end{gather*}
Furthemore, given $f\in\Hom_{\cal G}(V,W)$, let~$f_\Gamma:=f|_{V_\Gamma}$
and let $f^\Gamma$ be the induced map $V^\Gamma\to W^\Gamma$. It follows from the definitions that $f_\Gamma$ and
$f^\Gamma$ are morphisms of graded vector spaces and $\lie g$-modules.

\begin{prop} Let $V\in\Ob{\wh{\cal G}}$ and let $\Gamma$ be an interval closed subset of $\Lambda$.
\begin{enumerate}[{\rm(i)}]
\item\label{CAT175.i} Suppose that for each
$(\lambda,r)\in\Lambda(V)\setminus\Gamma$ there exists
$(\mu,s)\in\Gamma$ with $(\lambda,r)\prec(\mu,s)$.  Then $
V_\Gamma\in\Ob\wh{\cal G}[\Gamma]$ and $V/V_\Gamma\in\Ob\wh{\cal
G}[\Lambda\setminus\Gamma]$.

\item\label{CAT175.ii}Suppose that for each
$(\lambda,r)\in\Lambda(V)\setminus\Gamma$ there exists $(\mu,s)\in
\Gamma$ with $(\mu,s)\prec(\lambda,r)$. Then $
V_{\Lambda\setminus\Gamma}\in\Ob\wh{\cal
G}[\Lambda\setminus\Gamma]$ and  $ V^\Gamma\in\Ob\wh{\cal
G}[\Gamma]$.

\item\label{CAT175.iii} Let
$$
0\longrightarrow V\stackrel{f}{\longrightarrow} U\stackrel{g}\longrightarrow W\longrightarrow 0
$$
be a short exact sequence in~$\cal G$. Then~$U$ satisfies \eqref{CAT175.i} {\em(}respectively, \eqref{CAT175.ii}{\em)} if and only if~$V$
and $W$ satisfy \eqref{CAT175.i} {\em(}respectively, \eqref{CAT175.ii}{\em)} and
\begin{equation}\label{CAT175.10a}
0\longrightarrow V_\Gamma\stackrel{f_\Gamma}{\longrightarrow} U_\Gamma\stackrel{ g_\Gamma}{\longrightarrow}
W_\Gamma\longrightarrow 0
\end{equation}
{\em(}respectively, \begin{equation}\label{CAT175.10b} 0\longrightarrow V^\Gamma\stackrel{f^\Gamma}{\longrightarrow} U^\Gamma\stackrel{
g^\Gamma}{\longrightarrow} W^\Gamma\longrightarrow 0)
\end{equation} is an  exact sequence of objects in $\cal G[\Gamma]$.
\end{enumerate}
\end{prop}
\begin{pf} Consider the map $\lie g\otimes V_\Gamma\to V$ given by
$x\otimes v\mapsto (xt)v$. This is clearly a map of $\lie g$-modules.
Let $U$ be a $\lie g$-module complement to $V_\Gamma$ in $V$.
Suppose that  $V_\Gamma$ is not a subobject of $V$ in $\wh{\cal G}$, i.e
there exists $x\in\lie g$ and $v\in V_\Gamma$ such that
$(xt)v\notin V_\Gamma$. Since $(\lie g\otimes t)\bu(\lie
g)\subset\bu(\lie g)(\lie g\otimes t)$, we may
assume without loss of generality that $v\in V_\Gamma^+\cap
V[r]_\lambda$ for some $(\lambda,r)\in\Lambda$ and hence $\bu(\lie
g)v\cong_{\lie g} V(\lambda,r)$. In other words, the induced map of
$\lie g$-modules  $\lie g\otimes V(\lambda,r) \to U$ is
non-zero and so there exists $(\nu,r+1)\notin\Gamma$ such that the
composite map $\lie g\otimes V(\lambda,r) \to U\to
V(\nu,r+1)$ is non-zero.
This implies immediately that
$(\lambda,r)\prec(\nu,r+1)$. Choose $(\mu,s)\in\Gamma$ such that
$(\nu,r+1)\prec (\mu,s)$ for some $(\mu,s)\in\Gamma$. This gives~$(\nu,r+1)\in[(\lambda,r),(\mu,s)]$
which is   impossible since $\Gamma$ is
interval closed. Hence
$$(xt)v\in V_\Gamma,\qquad\forall \, x\in\lie g,\, v\in V_\Gamma^+,
$$ and \eqref{CAT175.i}  is proved.

In order to prove the second part, assume that~$V_{\Lambda\setminus\Gamma}$ is not a $\wh{\cal G}$-subobject of~$V$.
Then~$(xt)v\notin V_{\Lambda\setminus\Gamma}$ for some~$x\in\lie g$, $v\in V_{\Lambda\setminus\Gamma}$ and as before we may
assume, without loss of generality that~$v\in V_{\Lambda\setminus\Gamma}{}^+\cap V[r]_\lambda$ for some~$(\lambda,r)\in\Lambda\setminus
\Gamma$. Let~$U'$ be a $\lie g$-module complement of~$V_{\Lambda\setminus\Gamma}$. Then we have a non-zero $\lie g$-module
map $\lie g\tensor V(\lambda,r)\to U'$, given by extending $x\tensor v\mapsto (xt)v$, and hence a non-zero $\lie g$-module map
$\lie g\tensor V(\lambda,r)\to V(\nu,r+1)$ for some~$(\nu,r+1)\in \Gamma$. By assumption, there exists~$(\mu,s)\in\Gamma$ such
that $(\mu,s)\prec (\lambda,r)$. Thus, $(\lambda,r)\in [(\mu,s),(\nu,r+1)]\subset\Gamma$ since $\Gamma$ is interval closed,
which is a contradiction.

The first statement of~\eqref{CAT175.iii} is obvious since $\Lambda(U)=\Lambda(V)\cup\Lambda(W)$. For the second, note that we have
$$U\cong V\oplus W,$$ as graded $\lie g$-modules and hence as $\lie g$-module we have $$U_\Gamma\cong V_\Gamma\oplus W_\Gamma.$$ Since
$V_\Gamma$, $W_\Gamma$ and~$U_\Gamma$ are objects in $\cal G[\Gamma]$ the result follows.
\end{pf}
\begin{rem} Part \eqref{CAT175.i} of  the preceding proposition holds
for all $V\in\Ob\cal G$ and $\Gamma\subset\Lambda$ such that
$\Lambda(\soc(V))\subset\Gamma$ while part \eqref{CAT175.ii} holds if
$\Lambda(\head(V))\subset\Gamma$. In fact,  for any
$(\lambda,r)\in\Lambda(V)$, there exists
$(\mu,s)\in\Lambda(\soc(V))$ and $(\nu,p)\in\Lambda(\head(V))$ such
that
$$(\nu,p)\preccurlyeq(\lambda,r)\preccurlyeq(\mu,s).$$ This is an immediate
consequence of the fact that $V$ embeds in the injective envelope of
$\soc(V)$ and is a quotient of $P(\head(V))$ together with
Proposition~\ref{CAT130}\eqref{CAT130.iii} and Proposition \ref{injG}\eqref{injG.ii}.
Moreover, in this case if we let $\overline{V_\Gamma}$ be the
maximal subobject of $V$ that is in $\cal G[\Gamma]$, then the
Proposition implies that $$\overline{V_\Gamma}\cong V_\Gamma$$ if for
each $(\lambda,r)\in\Lambda(V)\setminus\Gamma$ there exists
$(\mu,s)\in\Gamma$ with $(\lambda,r)\preccurlyeq(\mu,s)$ and similarly
for $\overline V^\Gamma$.
\end{rem}

\subsection{}\label{projinjgamma} We isolate some consequences of the preceding Proposition
since we use them repeatedly in the following sections.
\begin{prop} Let $\Gamma$ be finite and interval closed and assume that
$(\lambda,r), (\mu,s)\in\Gamma$.
\begin{enumerate}[{\rm(i)}]
\item\label{projinjgamma.i} The object
$I(\lambda,r)_\Gamma$ is the injective envelope of $V(\lambda,r)$ in $\cal G[\Gamma]$ while $P(\lambda,r)^\Gamma$ is the projective cover
of~$V(\lambda,r)$ in~$\cal G[\Gamma]$.  In particular, $\cal G[\Gamma]$ has enough projectives.

\item\label{projinjgamma.ii}
We have
$$
[P(\lambda,r)^\Gamma:V(\mu,s)]=[P(\lambda,r):V(\mu,s)]=[I(\mu,s):V(\lambda,r)]=[I(\mu,s)_\Gamma:V(\lambda,r)].
$$
\item\label{projinjgamma.iii} For all $j\ge  0$, we have $$\Ext^j_{\cal G}(V(\lambda,r),
V(\mu,s))\cong\Ext^j_{\cal G[\Gamma]}(V(\lambda,r), V(\mu,s)).$$
\item\label{projinjgamma.iv}
Let $\boldsymbol p^\Gamma(\lambda,r):P(\lambda,r)\twoheadrightarrow
P(\lambda,r)^\Gamma$ {\em(}respectively, $\boldsymbol
\iota_\Gamma(\lambda,r):I(\lambda,r)_\Gamma \hookrightarrow
I(\lambda,r)${\em )} be the canonical projection {\em(}respectively,
the canonical embedding{\em)}. There exists an isomorphism
$\Hom_{\wh{\cal G}}(P(\lambda,r), P(\mu,s))\to\Hom_{\cal
G[\Gamma]}(P(\lambda,r)^\Gamma, P(\mu,s)^\Gamma)$ given by $f\to
f^\Gamma$ such that
$$\boldsymbol p^\Gamma(\mu,s)\circ f=f^\Gamma\circ\boldsymbol p^\Gamma(\lambda,r),$$ and
similarly an isomorphism $\Hom_{\cal G}(I(\mu,s),
I(\lambda,r))\to\Hom_{\cal G[\Gamma]}(I(\mu,s)_\Gamma,
I(\lambda,r)_\Gamma)$ given by $g\to g_\Gamma$ such that
\begin{equation*}
g\circ\boldsymbol \iota_\Gamma(\mu,s)=\boldsymbol \iota_\Gamma(\lambda,r)\circ g_\Gamma.\tag*{\qed}
\end{equation*}
\end{enumerate}
\end{prop}
\begin{pf}
It follows from~\propref{injG}\eqref{injG.iii} that
$$
(\nu,k)\prec (\lambda,r)\qquad \forall\,(\nu,k)\in\Lambda(I(\lambda,r))\backslash \{(\lambda,r)\}.
$$
\propref{CAT175}\eqref{CAT175.i} now gives,
$$I(\lambda,r)_\Gamma\in\Ob\cal G[\Gamma],\qquad
\soc(I(\lambda,r)_\Gamma)=V(\lambda,r)$$ 
and 
$$\Hom_{\cal G}(V(\mu,s),I(\lambda,r)/I(\lambda,r)_\Gamma)=0,\qquad \forall\,
(\mu,s)\in\Gamma.
$$ 
It follows immediately that~$\Ext^1_{\cal G[\Gamma]}(V(\mu,s),I(\lambda,r)_\Gamma)=0$ for all
$(\mu,s)\in\Gamma$ which implies the first statement
in~\eqref{projinjgamma.i}. The proof of the second statement is
similar.
The first and the last equality in~\eqref{projinjgamma.ii}
follow immediately from~\propref{CAT175} while the second equality
is an obvious consequence
of~\propref{CAT130}(\ref{CAT130.ii},\ref{CAT130.iv}) and
\propref{injG}.

Part~\eqref{projinjgamma.iii} is obvious if~$j=0$.
Set~$Q_{-1}(\mu,s)=V(\mu,s)$. For~$j\ge 0$ define inductively the
objects~$I_j(\mu,s)$ as the injective envelope of~$Q_{j-1}(\mu,s)$
in~$\cal G$ and~$Q_j(\mu,s)=\operatorname{coker}(Q_{j-1}(\mu,s)\hookrightarrow
I_j(\mu,s))$. Then
$$
0\to V(\mu,s)\to I_0(\mu,s)\to I_1(\mu,s)\to \cdots \to I_k(\mu,s)\to 0
$$
is an injective resolution for~$V(\mu,s)$ in~$\cal G$ and
\begin{equation}\label{projinjgamma.10}
\Ext^j_{\cal G}(V(\lambda,r),V(\mu,s))
\cong \Hom_{\cal G}(V(\lambda,r),I_j(\mu,s)),\qquad j>0.
\end{equation}
It follows from~\corref{injG.cor}\eqref{injG.cor.ii} by a straightforward induction on~$j$ that
$$(\nu,k)\not=(\mu,s)\in\Lambda(I_j(\mu,s))\cup \Lambda(Q_j(\mu,s))\implies
(\nu,k)\prec (\mu,s).
$$ 
Hence by \propref{CAT175},
$$Q_j(\mu,s)_\Gamma, I_j(\mu,s)_\Gamma\in\Ob\cal G[\Gamma], \qquad
I_j(\mu,s)/I_j(\mu,s)_\Gamma\in\cal G[\Lambda\setminus\Gamma]$$
and the sequence
\begin{equation}\label{projinjgamma.15}
0\to V(\mu,s)\to I_0(\mu,s)_\Gamma\to I_1(\mu,s)_\Gamma\to \cdots\to I_k(\mu,s)_\Gamma\to 0
\end{equation}
is exact in~$\cal G[\Gamma]$. It follows from part~\eqref{projinjgamma.i} that~\eqref{projinjgamma.15} is
an injective resolution of~$V(\mu,s)$ in~$\cal G[\Gamma]$.
Furthermore, for all~$(\lambda,r)\in\Gamma$
\begin{equation}\label{projinjgamma.20}
\Hom_{\cal G}(V(\lambda,r),I_j(\mu,s))\cong \Hom_{\cal G}(V(\lambda,r),I_j(\mu,s)_\Gamma)
\end{equation}
and similarly
$$
\Hom_{\cal G}(V(\lambda,r),Q_j(\mu,s))\cong \Hom_{\cal G}(V(\lambda,r),Q_j(\mu,s)_\Gamma).
$$
In particular, this implies that
$$
\soc(I_j(\mu,s)_\Gamma)\cong\soc(Q_{j-1}(\mu,s)_\Gamma)
$$
and so~$I_j(\mu,s)_\Gamma$ is the injective envelope
of~$Q_{j-1}(\mu,s)_\Gamma$ in~$\cal G[\Gamma]$.
By~\propref{CAT175}\eqref{CAT175.iii},
$I_j(\mu,s)_\Gamma/Q_{j-1}(\mu,s)_\Gamma\cong Q_j(\mu,s)_\Gamma$.
Then
$$\Ext^{j}_{\cal G[\Gamma]}(V(\lambda,r),V(\mu,s))
\cong \Hom_{\cal G[\Gamma]}(V(\lambda,r),I_j(\mu,s)_\Gamma)
$$
and~\eqref{projinjgamma.iii} follows from~\eqref{projinjgamma.10}
and~\eqref{projinjgamma.20}.

To prove \eqref{projinjgamma.iv},
let $f\in\Hom_{\wh{\cal G}}(P(\lambda,r), P(\mu,s))$, $f\not=0$. 
\iffalse
Since
$f(P(\lambda,r)_{\Lambda\setminus\Gamma})\subset
P(\mu,s)_{\Lambda\setminus\Gamma}$, we have an induced map
$f^\Gamma\in\Hom_{\cal G[\Gamma]}(P(\lambda,r)^\Gamma,
P(\mu,s)^\Gamma)$. Moreover
\fi
Then $f^\Gamma\ne 0$ since $f^\Gamma(1\otimes
V(\lambda,r))=f(1\otimes V(\lambda,r))\pmod{P(\lambda,r)_{\Lambda\setminus\Gamma}}\ne 0$. Thus, we have an injective
map
$$
\Hom_{\wh{\cal G}}(P(\lambda,r),P(\mu,s))\to \Hom_{\cal G[\Gamma]} (P(\lambda,r)^\Gamma,P(\mu,s)^\Gamma).
$$
Since both spaces have the same dimension by~\eqref{projinjgamma.ii} and~\propref{CAT150}, the isomorphism follows.
To prove the statement for injectives, observe that the natural map
$$\Hom_{\cal G}(I(\mu,s),I(\lambda,r))\to
\Hom_{\cal G}(I(\mu,s)_\Gamma,I(\lambda,r))$$
is surjective, while~$\Hom_{\cal G}(I(\mu,s)_\Gamma,I(\lambda,r)_\Gamma)
\cong \Hom_{\cal G}(I(\mu,s)_\Gamma,I(\lambda,r))$. The assertion follows from~\eqref{projinjgamma.ii}.
\end{pf}

\section{Algebras associated with the category~$\cal G$}\label{TMPALG}

In this section, we let $\Gamma$ be a finite interval closed
subset of~$\Lambda$. Given a finite dimensional algebra $A$, let
$A-\operatorname{mod}_f$ (respectively, $\operatorname{mod}_f-A$)
be the category of finite dimensional left (respectively, right)
$A$-modules.

\subsection{The algebra $\mathfrak A(\Gamma)$ and an equivalence of categories}\label{TMPALG10}
Set
$$
I(\Gamma)=\bigoplus_{(\lambda,r)\in\Gamma} I(\lambda,r),\qquad
\mathfrak A(\Gamma)=\End_{\cal G} I(\Gamma).
$$
Then $\mathfrak A(\Gamma)$ is an associative algebra. Moreover,
it is immediate from~\propref{projinjgamma} that
\begin{equation}\label{agamma}\mathfrak A(\Gamma)\cong \mathfrak A_\Gamma(\Gamma):=\End_{\cal G} I(\Gamma)_\Gamma.
\end{equation}
In particular, $\mathfrak A(\Gamma)-\operatorname{mod}_f$ is equivalent to $\mathfrak A_\Gamma(\Gamma)-\operatorname{mod}_f$
and similarly for the categories of right modules. Since $I(\Gamma)_\Gamma$ is the
injective cogenerator of $\cal G[\Gamma]$, a standard argument now
shows that the contravariant functor $\Hom_{\cal G}(-,I(\Gamma)_\Gamma)$ from~$\cal
G[\Gamma]$  to the category $\mathfrak A_\Gamma(\Gamma)-\operatorname{mod}_f$ is exact and provides a duality of categories.
Similarly, the functor $\Hom_{\cal G}(-,I(\Gamma)_\Gamma)^*$ from $\cal G[\Gamma]$ to the
category $\operatorname{mod}_f-\mathfrak A_\Gamma(\Gamma)$ is exact and provides an equivalence of categories.
Thus, $\cal G[\Gamma]$ is equivalent to~$\operatorname{mod}_f-\mathfrak A(\Gamma)$ and is dual to
$\mathfrak A(\Gamma)-\operatorname{mod}_f$.

It is clear from the definition that the simple objects in~$\mathfrak
A(\Gamma)-\operatorname{mod}_f$ are one-dimensional, that is to say
$\mathfrak A(\Gamma)$ is basic, and their isomorphism classes are
parametrized by elements of~$\Gamma$. Given~$(\lambda,r)\in\Gamma$,
let~$S_{\lambda,r}$ be the corresponding simple left $\mathfrak
A(\Gamma)$-module.
\begin{prop}\label{extalg}Let $\Gamma\subset\Lambda$ be finite and interval closed.
For all~$(\lambda,r),(\mu,s)\in\Gamma$, we have
\begin{equation}\label{sim}
\dim \Ext^1_{\mathfrak
A(\Gamma)}(S_{\lambda,r},S_{\mu,s})=\delta_{r,s+1}\dim \Hom_{\lie
g}(V(\mu),\lie g\tensor V(\lambda)).
\end{equation}
In particular,
the algebra $\frak A(\Gamma)$ is  quasi-hereditary.
\end{prop}
\begin{pf} Since $\Gamma$ is interval closed, it follows from \corref{projinjgamma}
that if $(\mu,s), (\lambda,r)\in\Gamma$, then  $$\Hom_{\cal G[\Gamma]}(V(\mu,s),
I(\lambda,r)_\Gamma/V(\lambda,r))\cong\Hom_{\cal G}(V(\mu,s),
I(\lambda,r)/V(\lambda,r)).$$ Hence, we have,
$$\Ext^1_{\mathfrak
A(\Gamma)}(S_{\lambda,r},S_{\mu,s})\cong\Ext^1_{\cal
G[\Gamma]}(V(\mu,s), V(\lambda,r))\cong\Ext^1_{\cal
G}(V(\mu,s),V(\lambda,r)).
$$
Equation~\eqref{sim}  follows from \propref{CAT170} which also proves that $\mathfrak
A(\Gamma)-\operatorname{mod}_f$ is a directed highest weight
category with poset of weights $(\Gamma,\preccurlyeq^{op})$.
Now \cite[Theorem~3.6]{CPS} implies that the algebra
$\mathfrak A(\Gamma)$ is quasi-hereditary.
\end{pf}

\subsection{}\label{TMPALG20}
Let~$Q(\Gamma)$ be the the $\Ext$-quiver of~$\mathfrak A(\Gamma)$,
that is, the quiver whose set of vertices is~$\Gamma$ and the number
of arrows from~$(\lambda,r)$ to~$(\mu,s)$ in $Q(\Gamma)$
is~$\dim\Ext^1_{\mathfrak A(\Gamma)}(S_{\lambda,r},S_{\mu,s})$. Note
that  that the number of paths from $(\lambda,r)$ to $(\mu,s)$ is
non-zero only if $(\mu,s)\prec(\lambda,r)$. In particular,
$Q(\Gamma)$ has no oriented loops. Let $\bc Q(\Gamma)$ be the path
algebra of $Q(\Gamma)$ and $\bc Q(\Gamma)[k]$ be the subspace
spanned by all paths of length $k$. Then
$$\bc Q(\Gamma)=\bigoplus_{k\in\bz_+}\bc Q(\Gamma)[k],$$
is a tightly graded associative algebra. Since~$\mathfrak A(\Gamma)$
is basic, a classical result of Gabriel's (cf. for
example~\cite[2.1(2)]{RinBook}) proves that  $\mathfrak
A(\Gamma)$ is isomorphic to a quotient of the path algebra~$\bc
Q(\Gamma)$ of~$Q(\Gamma)$ by an ideal~$R(\Gamma)$ which is contained
in the ideal of paths of length at least two. In particular this
means that an arrow between $(\lambda,r)$ and $(\mu,r-1)$ maps to a
non-zero element of $\Hom_{\cal G}(I(\lambda,r), I(\mu,r-1))$.
Given~$(\lambda,r)\in\Gamma$, let~$1_{\lambda,r}$ be the
corresponding primitive idempotent in~$\bc Q(\Gamma)$.
Note that $1_{\lambda,r}$ maps to the element
$\id_{\lambda,r}\in\End_{\cal G} I(\Gamma)$  defined by
$\id_{\lambda,r}(I(\mu,s))=\delta_{(\lambda,r), (\mu,s)}\id$ for
$(\mu,s)\in\Gamma$. In particular, $\id_{\mu,s} \mathfrak A(\Gamma)
\id_{\lambda,r}\cong \Hom_{\cal G}(I(\lambda,r),I(\mu,s))$ as a
vector space.

\subsection{A grading on $\mathfrak A(\Gamma)$}\label{TMPALG30}
Given $k\le r\in\bz_+$ define
\begin{gather*}
\mathfrak A(\Gamma)[k]=\bigoplus_{(\lambda,r),(\mu,r-k)\in\Gamma} \Hom_{\cal
G}(I(\lambda,r),I(\mu,r-k)).
\end{gather*}
Since $[I(\lambda,r): V(\mu,s)]=0$ unless $(\mu,s)\preccurlyeq(\lambda,r)$ (cf.~\propref{injG}\eqref{injG.iii}), it follows immediately that
$$
\mathfrak A(\Gamma)=\bigoplus_{k\in\bz_+} \mathfrak A(\Gamma)[k],\qquad \mathfrak A(\Gamma)[j] \mathfrak A(\Gamma)[k]\subset \mathfrak
A(\Gamma)[j+k],\quad\forall\, j,k\in\bz_+.
$$
Thus, $\mathfrak A(\Gamma)$ is a graded associative algebra and
$\mathfrak A(\Gamma)[0]$ is a commutative semi-simple subalgebra
of~$\mathfrak A(\Gamma)$. It is trivial to observe that with this
grading the algebra $\mathfrak A(\Gamma)$ is in fact a graded
quotient of $\bc Q(\Gamma)$ and hence  the ideal~$R(\Gamma)$ is
graded. In particular, $\mathfrak A(\Gamma)$ is tightly graded.

\subsection{The dimension of $R(\Gamma)$}\label{TMPALG40}
\begin{prop} Let $\Gamma$ be interval closed and finite and $(\lambda,r),(\mu,s)\in\Gamma$. Then the number of paths from $(\lambda,r)$ to $(\mu,s)$ is
$\dim\Hom_{\lie g}(V(\mu),\lie g^{\otimes (r-s)}\otimes V(\lambda))$.\end{prop}
\begin{pf}
Let $N(\lambda,r), (\mu,r-s))$ be the number of paths in $\bc
Q(\Gamma)$ from $(\lambda,r)$ to $(\mu,r-s)$ and
set~$N'(\lambda,\mu,s)=\dim\Hom_{\lie g}(V(\mu),\lie g^{\tensor
s}\tensor V(\lambda))$. It is easy to see that
$$
N'(\lambda,\mu,0)=\delta_{\lambda,\mu}=N((\lambda,r),(\mu,r)),
$$
while
$$
N'(\lambda,\mu,1)=\dim\Hom_{\lie g}(V(\mu),\lie g\tensor
V(\lambda))= \dim\Ext^1_{\mathfrak A(\Gamma)}
(S_{\lambda,r},S_{\mu,r-1})=N((\lambda,r),(\mu,r-1)),
$$
where we used~\eqref{sim}.
We now prove that ~$N$ and $N'$ satisfy the same recurrence relation
which establishes the proposition. It is clear that
$$
N((\lambda,r),(\mu,r-s))=\sum_{\nu\in P^+} N((\nu,r-s+1),(\mu,r-s))
N((\lambda,r),(\nu,r-s+1)).
$$
On the other hand, note that we can write
$$
\lie g^{\tensor (s-1)}\tensor V(\lambda)\cong\bigoplus_{\nu\in P^+}
V(\nu)^{N'(\lambda,\nu,s-1)}.
$$
Tensoring with $\lie g$ gives,
\begin{equation*}
N'(\lambda,\mu,s)=\dim\bigoplus_{\nu\in P^+}\Hom_{\lie g}(V(\mu),
\lie g\tensor V(\nu))^{N'(\lambda,\nu,s-1)}
\\=\sum_{\nu\in P^+} N'(\lambda,\nu,s-1) N'(\nu,\mu,1),
\end{equation*}
and the proof is complete.
\end{pf}

\begin{cor}\label{TMPALG40.cor}
 Given  $(\mu,r-s),
(\lambda,r)\in\Gamma$, we have
\begin{equation*}
\dim 1_{\mu,r-s}R(\Gamma) 1_{\lambda,r} =\dim\Hom_{\lie g}(V(\mu),\lie g^{\tensor s}\tensor V(\lambda))-\dim\Hom_{\lie g}(V(\mu),S^{(s)}(\lie
g)\tensor V(\lambda)).
\end{equation*}
In particular, the algebra $\mathfrak A(\Gamma)$ is hereditary if and only if
\begin{equation*} \dim\Hom_{\lie g}(V(\mu),\lie g^{\tensor s}\tensor
V(\lambda))=\dim\Hom_{\lie g}(V(\mu),S^{(s)}(\lie g)\tensor V(\lambda))
\end{equation*}
for all $(\mu,r-s),(\lambda,r)\in\Gamma$.
\end{cor}
\begin{pf} Observe that
\begin{equation*}
\begin{split}\dim 1_{\mu,r-s}R(\Gamma)
1_{\lambda,r}&=N((\lambda,r),(\mu,r-s))-\dim \id_{\mu,r-s}\mathfrak
A(\Gamma) \id_{\lambda,r}\\
&=N((\lambda,r),(\mu,r-s))-\dim\Hom_{\cal G}(I(\lambda,r),I(\mu,r-s)).
\end{split}
\end{equation*}
The first assertion is now immediate from the above Proposition and~\propref{injG}\eqref{injG.ii}. For the second,
it is enough to observe that $\mathfrak A(\Gamma)$ is hereditary if and only if $R(\Gamma)=0$ and that
$R(\Gamma)=\bigoplus_{(\lambda,r),(\mu,r-s)\in\Gamma} 1_{\mu,r-s} R(\Gamma) 1_{\lambda,r}$.
\end{pf}

\subsection{}\label{TMPALG50}
Given $\Gamma\subset \Lambda_{\le r}$, let~$\Gamma^{\#_r}=\{ (-w_\circ \mu,r-s)\,:\, (\mu,s)\in\Gamma\}$.
It is easy to see that $\Gamma^{\#_r}$ is interval closed if and only if $\Gamma$ is interval closed.
\begin{prop}
Suppose that $\Gamma\subset \Lambda_{\le r}$ is finite and interval closed. Then $\mathfrak A(\Gamma^{\#_r})
\cong \mathfrak A(\Gamma)^{op}$.
\end{prop}
\begin{pf}
Let~$(\lambda,r)\in\Gamma$. Then~$V(\lambda,r)^{\#_r}$ is an object in~$\cal G[\Gamma^{\#_r}]$ and it follows
that $(\cal G[\Gamma])^{\#_r}=\cal G[\Gamma^{\#_r}]$. Thus, $\cal G[\Gamma]$ is dual to $\cal G[\Gamma^{\#_r}]$.
It follows from~\ref{TMPALG10} that $\mathfrak A(\Gamma)^{op}$ and $\mathfrak A(\Gamma^{\#_r})$ are Morita
equivalent. Since they are both basic, the assertion follows.
\end{pf}

\section{Examples of $\Gamma$ with $\mathfrak A(\Gamma)$ hereditary}\label{EXH}

Throughout this section we use the notations of~\cite{RinBook} for
the various types of quivers. This should eliminate confusion with
the notation for the types of simple Lie algebras. For instance,
$\mathbb T_{n_1,\dots,n_r}$ denotes a quiver whose underlying graph
is a star with $r$ branches where the $i$th branch contains~$n_i$
vertices, while $\tilde{\mathbb X}_k$ denotes the quiver whose
underlying graph is of affine Dynkin type $X_k^{(1)}$. Notice that if~$\Gamma\subset\Lambda_{\le r}$, then~\propref{TMPALG50} implies
that $Q(\Gamma^{\#_r})$ is the opposite quiver of~$Q(\Gamma)$.

\subsection{The generalized Kronecker quivers} Let~$\lambda\in P^+$ be non-zero and
let $$k_\lambda=|\{ i\in I\,:\, \lambda(h_i)>0\}|.$$ It is easily
checked, by using Lemma \ref{CAT10} that
$$
\dim\Hom_{\lie g}(V(\lambda),\lie g\tensor
V(\lambda))=k_\lambda.$$For $r\in\bz_+$, set $\Gamma_{\lambda,r}=\{
(\lambda,r),(\lambda,r+1)\}$. It follows that
$Q(\Gamma_{\lambda,r})$ is the quiver with $k_\lambda$ arrows from
$(\lambda,r+1)$ to $(\lambda,r)$. If $k_\lambda=1$, the quiver is of
type $\mathbb A_2$, if $k_\lambda=2$ it is the Kronecker quiver
$\tilde{\mathbb A}_1$, while for $k_\lambda>2$ we get the
generalized Kronecker quiver. Since there are no paths of length two
in these quivers, it follows that $\mathfrak
A(\Gamma_{\lambda,r})\cong \bc Q(\Gamma_{\lambda,r})$.

\subsection{Quivers of type $\tilde{\mathbb D}_4$} Suppose that $\lie g$ is not of type~$\lie{sl}_2$.
Let~$I_{\bullet}=\{i\in I\,:\, \theta-\alpha_i\in R^+\}$. Note that $|I_{\bullet}|=1$
if $\lie g$ is not of type~$\lie{sl}_{n+1}$ and let~$i_\bullet$ be the unique element of~$I_{\bullet}$.
If $\lie g\cong\lie{sl}_{n+1}$, $n>1$, $I_\bullet=\{1,n\}$.

Let~$r\in\bz_+$. If $\lie g$ is not of type~$\lie{sl}_{n+1}$ set
$$
\Gamma=\{(\theta,r),(0,r+1), (2\theta,r+1),(2\theta-\alpha_{i_\bullet},r+1),
(\theta,r+1) \}.
$$
Otherwise, set
$$
\Gamma=\{(\theta,r), (0,r+1), (2\theta,r+1), (2\theta-\alpha_1,r+1),
(2\theta-\alpha_n,r+1)\}.
$$
In the first case, we find that  $Q(\Gamma)$ is
$$
\makeatletter
\def\dggeometry{\unitlength=0.018pt\dg@YGRID=2\dg@XGRID=3}
\def\dgeverynode{\scriptstyle}
\begin{diagram}
\node{}\node{(2\theta-\alpha_{i_\bullet},r+1)}\arrow{s}\\
\node{(0,r+1)}\arrow{e}\node{(\theta,r)}\node{(2\theta,r+1)}\arrow{w} \\
\node{}\node{(\theta,r+1)}\arrow{n}
\end{diagram}
$$
while in the second case $Q(\Gamma)$ is
$$ \makeatletter
\def\dggeometry{\unitlength=0.018pt\dg@YGRID=2\dg@XGRID=3}
\def\dgeverynode{\scriptstyle}
\begin{diagram}
\node{}\node{(2\theta-\alpha_{1},r+1)}\arrow{s}\\
\node{(0,r+1)}\arrow{e}\node{(\theta,r)}\node{(2\theta,r+1)}\arrow{w} \\
\node{}\node{(2\theta-\alpha_n,r+1)}\arrow{n}
\end{diagram}
$$
 The algebra $\mathfrak A(\Gamma)$ is hereditary
since any path in $Q(\Gamma)$ has length at most one.
Note that $\Gamma$ can be shifted by any $\lambda\in P^+$ such that
$\dim\Hom_{\lie g}(V(\lambda+\theta),\lie g\tensor
V(\lambda+\theta))=1$.

\subsection{Quivers of type $\mathbb A_\ell$} If  $\lie g$ is not of type $\lie{sl}_2$ choose $i_\bullet\in I_\bullet$.

\begin{prop}\label{aquiver} Fix~$\lambda\in P^+$ with $\lambda(h_{i_\bullet})\ne 0$, $\ell\in\bz_+$. Let $r_j\in\bz_+$, $0\le j\le \ell$ be
such that $|r_k-r_{k+1}|=1$ for $0\le k\le \ell-1$. Let
$\alpha\in\{\theta, \theta-\alpha_{i_\bullet}\}$. The set
$$
\Gamma=\{(\lambda+j\alpha, r_j):0\le j\le \ell
\}
$$
is interval closed, the quiver $Q(\Gamma)$ is of type
 $\mathbb A_{\ell+1}$ and the algebra $\mathfrak A(\Gamma)$ is
 hereditary.
\end{prop}
 \begin{pf} We prove the proposition in the case when $\alpha=
\theta-\alpha_{i_\bullet}$, the proof in the other case being similar
and in fact simpler. Suppose that for some $0\le j,j'\le\ell$ and
$(\mu,s)\in\Lambda$, we have,
$$(\lambda+j\alpha,r_j)\prec(\mu,s)\prec
(\lambda+j'\alpha,r_{j'}),$$ then we have
$$\mu=\lambda+j'\alpha-\sum_{p=1}^{r_{j'}-s}\beta_p=\lambda+j\alpha+\sum_{q=1}^{s-r_j}\gamma_q,$$ and $$ r_j<s<r_{j'}.$$
for some $\beta_p, \gamma_q\in R\sqcup\{0\}$. Assuming without loss
of generality that $j\le j'$, we find in particular,
$0<r_{j'}-r_j\le j'-j$ and
$$ (j'-j)(\theta-\alpha_{i_\bullet})=\sum_{p=1}^{r_{j'}-s}
\beta_p+\sum_{q=1}^{s-r_j}\gamma_q.
$$  Equating the coefficients of
$\alpha_{i}$, $i\not=i_\bullet$ on both sides of  the above
expression we conclude that $$\beta_p,\gamma_q\in
\{\theta,\theta-\alpha_{i_\bullet}\},  \ \ \ r_{j'}-r_j=j'-j,$$ for
$1\le p\le r_{j'}-s$ and $1\le q\le s-r_j$, which gives
$$r_p=r_j+(p-j),\ \ j\le p\le j'.$$
 Next, equating the
coefficient of~$\alpha_{i_\bullet}$  on both sides now gives
$\beta_p=\gamma_q=\theta-\alpha_{i_\bullet}$ for all $1\le p\le
r_{j'}-s$, $1\le q\le s-r_j$. This proves that
$$(\mu,s)=(\lambda+s\alpha,s)=(\lambda+s\alpha, r_s),$$ and hence
$(\mu,s)\in\Gamma$.

It follows from~\propref{extalg} that
$$\dim\Ext^1_{\mathfrak A(\Gamma)}(S_{\lambda+j\alpha, r_j}, S_{\lambda+ k\alpha,
r_{k}})=\delta_{r_{j}-r_k,1}\dim\Hom_{\lie
g}(V(\lambda+k\alpha),\lie g\otimes V(\lambda+j\alpha)), $$ and
applying Lemma \ref{CAT10} now gives
$$\dim\Ext^1_{\mathfrak A(\Gamma)}(S_{\lambda+j\alpha, r_j}, S_{\lambda+ k\alpha,
r_{k}})=\delta_{r_{j}-r_k,1}\delta_{|j-k|,1}.$$
 This shows that there is
precisely one arrow between $(\lambda+j\alpha,r_j)$ and
$(\lambda+(j\pm 1)\alpha,r_{j\pm1})$ and no other arrow which has
$(\lambda+j\alpha,r_j)$ as its head or tail. Therefore,
$Q(\Gamma)$ is of type~$\mathbb A_{\ell+1}$.

To prove that the algebra $\mathfrak A(\Gamma)$ is hereditary, let
$(\lambda+k\alpha,r_k), (\lambda+k'\alpha, r_{k'})\in \Gamma$. The
number of paths in $Q(\Gamma)$ between these vertices is zero unless
$(\lambda+k\alpha, r_k)$, $(\lambda+k'\alpha, r_{k'})$ are strictly
comparable.  Assume without loss of generality that
$(\lambda+k\alpha,r_k)\prec (\lambda+k'\alpha,r_{k'})$ and also that
$k\le k'$. But in this case, we have proved that $r_k=r_{k'}-k'+k$
and that there is exactly one path from $(\lambda+k'\alpha,r_{k'})$
to $(\lambda+k\alpha,r_{k'}-k'+k)$. The result now follows from
\corref{TMPALG40} if we prove that
$$\dim\Hom_{\lie g}(V(\lambda+k\alpha), S^{(k'-k)}(\lie g)\otimes
V(\lambda+k'\alpha))=1.
$$
Since $\Hom_{\lie g}(V((k-k')\theta)),S^{(k'-k)}(\lie g))\ne 0$ it
suffices to prove that $$\dim\Hom_{\lie g}(V(\lambda+k\alpha),
V((k-k')\theta))\otimes V(\lambda+k'\alpha))=1.
$$
But this again follows from Lemma~\ref{CAT10}.
\end{pf}
\begin{rem} The restriction $\lambda(h_{i_\bullet})\ne 0$ is  not necessary if
$\alpha=\theta$.\end{rem}

\subsection{Quivers of
type $\tilde{\mathbb D}_{\ell+1}$, $\ell\ge 4$} The arguments given
in the previous section can be used with obvious modifications to
prove the following.
 Let $\alpha=\theta-\alpha_{i_\bullet}$, $\ell\in\bz_+$ with
$\ell\ge 4$ and $\lambda\in P^+$.
 Let $r_j\in\bz_+$, $2\le j\le \ell-1$ be
such that $|r_k-r_{k+1}|=1$ for $2\le k\le \ell-2$ and let
$$\Gamma=\Gamma_1\cup\Gamma_2$$ where $\Gamma_1= \{(\lambda+j\theta, r_j):2\le j\le \ell-1
 \}$ and $$\Gamma_2=\{(\lambda+\theta, r_2+1),\, (\lambda+\theta+\alpha, r_2+1),\, (\lambda+\ell\theta, r_{\ell-1}-1),\,
( \lambda+(\ell-1)\theta+\alpha, r_{\ell-1}-1)\}.
$$
Then $\Gamma$ is interval closed and the quiver $Q(\Gamma)$ is
of type $\tilde{\mathbb D}_{\ell+1}$, where the set $\Gamma_2$ consists of precisely
 those vertices which are either the head or the tail of precisely
 one arrow. Moreover, the algebra $\mathfrak A(\Gamma)$ is
 hereditary.

\subsection{Star-shaped quivers} Suppose that $\lie g$ is not of type~$C_n$.
Fix~$\lambda\in P^+$. For $1\le p\le 3$ and $0\le j\le \ell_p$,
let~$r_{p,j}\in\bz_+$ be such that  $|r_{p,j}-r_{p,j+1}|=1$. Let
\begin{align*}
\Gamma_1&=\{ (\lambda+(\ell_1-j)\theta,r_{1,j}): 0\le j\le \ell_1\},\\
\Gamma_2&=\{ (\lambda+(\ell_1+j)\theta,r_{2,j}): 0\le j\le \ell_2\},\\
\Gamma_3&=\{
(\lambda+(\ell_1+j)\theta-j\alpha_{i_\bullet},r_{3,j}):0\le j\le
\ell_3-1\}.
\end{align*}
Set~$\Gamma=\Gamma_1\cup\Gamma_2\cup\Gamma_3$.
\begin{prop}
Suppose that for all $1\le p\le 3$ and~$0\le j_p<\ell_p$ we have
$$|r_{3,j_3}-r_{2,j_2}|\le |j_3-j_2|,\ \  |r_{3,j_3}-r_{1,j_1}|\le
|j_3-j_1|.$$  Then~$\Gamma$ is interval closed, $Q(\Gamma)$ is of
type~$\mathbb T_{\ell_1,\ell_2,\ell_3}$ and~$\mathfrak A(\Gamma)$ is
hereditary.
\end{prop}
\begin{pf}
Note that the conditions  imply that
$$r_{1,j}=r_{2,j}=r_{3,j},\quad 0\le j\le \ell_3,\qquad
\Gamma_1\cap\Gamma_2\cap\Gamma_3=\{(\lambda+\ell_1\theta,r_{1,0})\}.$$
Using  \propref{aquiver} we see that $\Gamma_1\cup\Gamma_2$ and
$\Gamma_3$ are  interval closed and that $\frak
A(\Gamma_1\cup\Gamma_2)$ and $\frak A(\Gamma_3)$ are hereditary of
type $\mathbb A_{\ell_1+\ell_2-1}$ and $\mathbb A_{\ell_3}$.
 The
proposition follows at once if we prove that an element in
~$\Gamma_3$ is not comparable in the partial order $\prec$ to any
element in~$(\Gamma_1\cup\Gamma_2)$ except possibly to
$(\lambda+\ell_1\theta,r_{1,0})$.

Suppose first that
$(\lambda+(\ell_1+j_3)\theta-j_3\alpha_{i_\bullet}, r_{3,j_3})$ is
strictly comparable with $(\lambda+(\ell_1+j_2)\theta,r_{2,j_2})$.
Then~$0<|r_{2,j_2}-r_{3,j_3}|\le
 |j_3-j_2|$ and
\begin{equation}\label{TMPPP.20}
(j_2-j_3)\theta+j_3\alpha_{i_\bullet}=\sum_{p=1}^{|r_{2,j_2}-r_{3,j_3}|}
\beta_p,\qquad \beta_p\in R\sqcup\{0\}.
\end{equation}
If~$j_3>j_2$, then we see by comparing the coefficients of
$\alpha_i$ with $i\ne i_\bullet$ on both sides, that $\beta_p\in
\{-\theta,-\theta+\alpha_{i_\bullet}\}$ and that
$|r_{2,j_2}-r_{3,j_3}|=j_3-j_2$. Suppose that $-\theta$ occurs $s$
times in the set $\{\beta_p:1\le p\le j_3-j_2\}$.  Then $0\le s\le
j_3-j_2$ and $-\theta+\alpha_{i_\bullet}$ occurs $j_3-j_2-s$ times.
It follows that $-(j_3-j_2)\theta+j_3\alpha_{i_\bullet}=\sum_p
\beta_p = -s\theta+(j_3-j_2-s)(-\theta+\alpha_{i_\bullet})
=-(j_3-j_2)\theta+(j_3-j_2-s)\alpha_{i_\bullet}$, which implies
$j_2=s=0$.

Furthermore, suppose that~$j_2>j_3$. By comparing the coefficients
of $\alpha_i$ with $i\ne i_\bullet$ in both sides
of~\eqref{TMPPP.20}, we conclude
that~$\beta_p\in\{\theta,\alpha_{i_\bullet}\}$ and
$|r_{2,j_2}-r_{3,j_3}|=j_2-j_3$. Suppose that $\alpha_{i_\bullet}$
occurs $s'$ times. Then we must have $s'=j_3$ and
$j_2-j_3-s'=(j_2-j_3)$ which is only possible if~$s'=j_3=0$.

Suppose now that
$(\lambda+(k+j_3)\theta-j_3\alpha_{i_\bullet},r_{3,j_3})$ is
strictly comparable with $(\lambda+{(k-j_1)}\theta,r_{1,j_1})$. Then
we must have $0\le r=|r_{3,j_3}-r_{1,j_1}|\le |j_3-j_1|$ and
$$
(j_1+j_3)\theta-j_3\alpha_{i_\bullet}=\sum_{p=1}^{r} \gamma_p,\qquad
\gamma_p\in R\sqcup\{0\}.
$$
Comparing the coefficients of~$\alpha_i$, $i\not=i_\bullet$ in both
sides shows that $\beta_p\in\{\theta,\theta-\alpha_{i_\bullet}\}$.
If $\theta$ appears $s$ times, then we have
$$
r\theta-(r-s)\alpha_{i_\bullet}=(j_1+j_3)\theta-j_3\alpha_{i_\bullet},
$$
which implies $r=j_1+j_3$ and~$r-s=j_3$. Since $r\le |j_1-j_3|$ the
first equality implies that either~$j_1=0$ or~$j_3=0$.
\end{pf}

\section{Quivers with relations}\label{QR}

In this section we give an example of $\Gamma$ for which $\frak
A(\Gamma)$ is not hereditary. The example is motivated by a family
of $\lie g[t]$-modules called the Kirillov-Reshetikhin modules
(cf.~\cite{CMkir1}).  We assume in this section that $\lie g$ is of
type $D_n$, $n\ge 6$. Recall that for $i\in I$ with $i\ne n-1,n$ we
have $$\omega_i=\sum_{j=1}^i j\alpha_j+
i\sum_{j=i+1}^{n-2}\alpha_j+\frac i2(\alpha_{n-1}+\alpha_n).$$
Let~$\Gamma$ be the interval $[(2\omega_4,0),(0,4)]$. It is easily
checked that $$\Gamma=\{(2\omega_4,0),(\omega_2+\omega_4,1),
(\omega_4,2),(2\omega_2,2),(\omega_1+\omega_3,2),(\omega_2,3),(0,4)\}.$$
Since $\lie g\cong_{\lie g} V(\omega_2)$, it is now not hard to see
by using \propref{sim} that the quiver $Q(\Gamma)$ is as follows:
$$
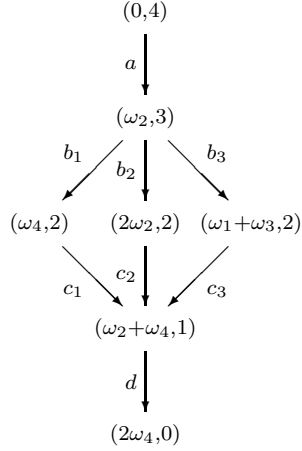
\makeatletter
\def\dggeometry{\unitlength=0.04pt\dg@YGRID=1\dg@XGRID=1}
\def\dgeverynode{\scriptstyle}
\begin{diagram}
\node{}\node{(0,4)}\arrow{s,l}{a}\\
\node{}\node{(\omega_2,3)}\arrow{se,l}{b_3}\arrow{s,l}{b_2}\arrow{sw,l}{b_1}\\
\node{(\omega_4,2)}\arrow{se,r}{c_1}\node{(2\omega_2,2)}\arrow{s,l}{c_2}\node{(\omega_1+\omega_3,2)}\arrow{sw,r}{c_3}\\
\node{}\node{(\omega_2+\omega_4,1)}\arrow{s,l}{d}\\
\node{}\node{(2\omega_4,0)}
\end{diagram}
$$
The path algebra $\bc
Q(\Gamma)$ has a basis consisting of the paths of length at most four
which we list below for the readers convenience:
\begin{alignat*}{3}&\{1_{\lambda,r}:(\lambda,r)\in\Gamma\},&\qquad
&\{a,b_i,c_i,d: 1\le i\le 3\},&\\
&\{b_ia,dc_i, c_ib_i:1\le i\le 3\},& &\{c_ib_ia, dc_ib_i:1\le i\le
3\},& \qquad \{dc_ib_ia: 1\le i\le 3\}
\end{alignat*}
We now compute the dimension of $\dim 1_{\mu,s} R(\Gamma)
1_{\lambda,r}$ for $(\mu,s),(\lambda,r)\in\Gamma$ with $r-s\ge 2$.
By \corref{TMPALG40} it suffices to calculate $\dim \Hom_{\cal
G}(I(\lambda,r),I(\mu,s))$.
Using~\propref{projinjgamma}\eqref{projinjgamma.ii} and the graded
characters of injective envelopes of simples in $\cal G[\Gamma]$
listed in Appendix~\ref{app.1} we find that
$$\dim 1_{\mu,s} R(\Gamma) 1_{\lambda,r}= 1,$$ if
$$
((\lambda,r),(\mu,s))\in\{((0,4),(\omega_1+\omega_3,2)),
((\omega_2,3),(\omega_2+\omega_4,1)),((\omega_1+\omega_3,2),(2\omega_4,0))\},$$
and $$\dim 1_{\mu,s} R(\Gamma) 1_{\lambda,r}= 2,$$ if
$$((\lambda,r),(\mu,s))\in\{((0,4),(\omega_2+\omega_4,1)),
((0,4),(2\omega_4,0))\},$$ while $\dim 1_{\mu,s} R(\Gamma)
1_{\lambda,r}= 0$ otherwise. This implies that there exists a unique
(up to multiplication by  non-zero constants) choice of  complex
numbers $x_i$, $1\le i\le 3$ and $\xi_j$, $\zeta_j$, $\eta_j$,
$j=1,2$ such that $\mathfrak A(\Gamma)$ is the quotient of $\bc
Q(\Gamma)$ by the following relations
\begin{equation}\label{quad}
b_3 a=0=d c_3,\qquad x_1 c_1 b_1+x_2 c_2 b_2+x_3 c_3 b_3 = 0,
\end{equation}
\begin{equation}\label{cub}
c_3 b_3 a =0,\qquad \xi_1 c_1 b_1 a +\xi_2 c_2 b_2 a=0,\qquad
\zeta_1 d c_1 b_1+\zeta_2 d c_2 b_2 =0,
\end{equation}
and \begin{equation}\label{quart}\eta_1 dc_1 b_1a+\eta_2 d c_2 b_2 a=0,\qquad
dc_3 b_3 a=0.
\end{equation}

\begin{prop} The algebra $\frak A(\Gamma)$ is the quotient of
$\bc Q(\Gamma)$ by the relations:
\begin{equation}\label{min} b_3a=0,\qquad dc_3=0,\qquad
c_1b_1+c_2b_2+c_3b_3=0.
\end{equation}
In particular $\frak
A(\Gamma)$ is quadratic, of global dimension $2$ and of tame
representation type.
\end{prop}
\begin{rem} It is not hard to see by using the results of \cite{CMkir1} and the equivalence of categories between $\mathfrak A(\Gamma)$-modules
and $\cal G[\Gamma]$, that the projective cover in $\mathfrak A(\Gamma)-\operatorname{mod}_f$ of $S_{0,4}$
or the injective envelope of~$S_{2\omega_4,0}$ corresponds to  the
Kirillov-Reshetikhin module $KR(2\omega_4)$, which is thus injective and projective in~$\cal G[\Gamma]$.
This connection will be
pursued elsewhere.
\end{rem}
\begin{pf} The relations in \eqref{min} are clearly independent of each other. To see that all relations in
$\mathfrak A(\Gamma)$ are consequences of those in~\eqref{min} it is enough to prove that the space spanned by
$b_ic_i$, $b_jc_j$ with $1\le i< j\le 3$ is always of dimension two.
Using the equivalence of categories, \propref{TMPALG50},
\propref{injG} and~\propref{CAT90}\eqref{CAT90.i} this can be
reformulated into the following question on morphisms in $\cal G$.
Thus, for $\mu\in\{2\omega_2, \omega_4, \omega_1+\omega_3\}$  fix
non-zero elements $f_\mu\in\Hom_{\wh{\cal
G}}(P(\omega_2+\omega_4,2), P(\mu,1))$ and $g_\mu\in\Hom_{\wh{\cal
G}}( P(\mu,1),P(\omega_2,0))$. We have to prove that the elements
$g_\mu f_\mu$ and $g_\lambda f_\lambda$ are linearly independent in
$\Hom_{\wh{\cal G}}(P(\omega_2+\omega_4,2),P(\omega_2,0))$. In turn,
using~\propref{CAT150} this question translates into the following
question in the category $\cal F(\lie g)$. Let $\bar f_\mu$, $\bar
g_\mu$ be the restrictions of $f_\mu$ and~$g_\mu$ to
$V(\omega_2+\omega_4)$ and $V(\mu)$ respectively and $\bop: T^3(\lie
g)\to S^2(\lie g)\otimes \lie g$ be the canonical projection. The
elements $\bop\circ(1\otimes\bar g_\mu)\circ \bar f_\mu$ and
$\bop\circ(1\otimes\bar g_\lambda)\circ \bar f_\lambda$ are linearly
independent elements of $\Hom_{\lie g}(V(\omega_2+\omega_4),
S^2(\lie g)\otimes \lie g)$. This is done by an explicit computation
of the maps, and the details can be found in the Appendix~\ref{app.2}.

Since $\mathfrak A(\Gamma)$ is quadratic, it follows from~\cite[Theorem~1.1]{Bon}
that $\Ext^2_{\mathfrak A(\Gamma)}(S_{\lambda,r},S_{\mu,s})=0$ unless~$r=s+2$. We have
$$
\dim\Ext^2_{\mathfrak A(\Gamma)}(S_{0,4},S_{\omega_1+\omega_3,2})=
\dim \Ext^2_{\mathfrak A(\Gamma)}(S_{\omega_2,3}, S_{\omega_2+\omega_4,1})=\\
\dim\Ext^2_{\mathfrak
A(\Gamma)}(S_{\omega_1+\omega_3,2},S_{2\omega_4,0})=1
$$
and $\Ext^2_{\mathfrak A(\Gamma)}(S_{\lambda,r},S_{\mu,r-2})=0$ in
all other cases. We claim that~$\Ext^3_{\mathfrak A(\Gamma)}(S_{\lambda,r},S_{\mu,s})=0$
for all~$(\lambda,r),(\mu,s)\in\Gamma$. Indeed, by~\cite[Theorem~1.1]{Bon}
\begin{multline*}
\dim\Ext^3_{\mathfrak A(\Gamma)}(S_{\lambda,r},S_{\mu,s})\\= \dim
1_{\mu,s} ((\bc Q(\Gamma)_+ R(\Gamma)\cap R(\Gamma)\bc
Q(\Gamma)_+)/(R(\Gamma)^2+\bc Q(\Gamma)_+ R(\Gamma)\bc Q(\Gamma)_+)))
1_{\lambda,r},
\end{multline*}
where~$\bc Q(\Gamma)_+$ is the radical of~$\bc Q(\Gamma)$.
If $r-s<4$, it is clear that $\dim 1_{\mu,s} (\bc Q(\Gamma)_+
R(\Gamma)\cap R(\Gamma)\bc Q(\Gamma)_+) 1_{\lambda,r}=0$. For
$r-s=4$, we have a unique pair $(\lambda,r),(\mu,r-4)\in\Gamma$,
namely $(0,4)$ and~$(2\omega_4,0)$, and two linearly independent
elements in $1_{2\omega_4,0}(\bc Q(\Gamma)_+ R(\Gamma) \cap R(\Gamma) \bc
Q(\Gamma)_+)1_{0,4}$, namely $d c_3 b_3 a$ and $d c_2 b_2 a+d c_1 b_1 a$.
The first is contained in~$R(\Gamma)^2$, since it can be written as
$(d c_3)(b_3 a)$ and $d c_3, b_3 a\in R(\Gamma)$, while the second
is contained in $\bc Q(\Gamma)_+ R(\Gamma)\bc Q(\Gamma)_+$ since it
can be written as $d(c_1 b_1+c_2 b_2+c_3 b_3)a$. Thus,
$\dim\Ext^3_{\mathfrak A(\Gamma)}(S_{0,4},S_{2\omega_4,0})=0$.

It remains to prove that the algebra is tame.
Let~$\Gamma_0=\Gamma\setminus\{(2\omega_4,0), (0,4)\}$. Note that~$\Gamma_0$
is interval closed and
consider the subalgebra $\mathfrak A(\Gamma_0)$. This
algebra is canonical (cf.~\cite[3.7]{Rin}) and of type~$(2,2,2)$,
hence tame concealed (\cite[4.3(5)]{Rin}). Let $K$ be the subspace
of~$\mathfrak A(\Gamma_0)$ spanned by~$\{b_3, c_3 b_3\}$. Clearly,
$K$ is a $\mathfrak A(\Gamma_0)$-submodule of~$\mathfrak
A(\Gamma_0)1_{\omega_2,3}$. Let $M$ be the quotient of $\mathfrak
A(\Gamma_0)1_{\omega_2,3}$ by $K$. This $\mathfrak
A(\Gamma_0)$-module has dimension vector
$$
\makeatletter
\def\dggeometry{\unitlength=0.025pt\dg@YGRID=1\dg@XGRID=1}
\def\dgeverynode{\scriptstyle}
\begin{diagram}
 \node{}\node{1}\arrow{se,-}{}\arrow{s,-}{}\arrow{sw,-}{}\\
\node{1}\arrow{se,-}{}\node{1}\arrow{s,-}{}\node{0}\arrow{sw,-}{}\\
\node{}\node{1}
\end{diagram}
$$
and hence belongs to the tubular family of type~$(2,2,2)$.
Then it is easy to check that $\mathfrak A(\Gamma)$ is obtained as
one-point extension and one-point coextension of $\mathfrak
A(\Gamma_0)$ at~$M$ and hence is tame (even domestic). We refer the reader to~\cite[4.7]{RinBook}
for details.
\end{pf}

\appendix

\section{}
\allowdisplaybreaks[4]
\subsection{}\label{app.1}
The graded characters of injectives envelopes of simple objects
in~$\cal G[\Gamma]$ given by~\propref{injG}
and~\propref{projinjgamma} can be calculated explicitly by using the
{\sf LiE} computer program (\cite{LiE}) and we list them below for
the reader's convenience.
\begin{align*}
I(0,4)_\Gamma&=V(0,4)\oplus V(\omega_2,3)\oplus \left(\begin{array}{c}V(\omega_4,2)\\ \oplus\\ V(2\omega_2,2)\end{array}\right)
\oplus V(\omega_2+\omega_4,1)\oplus V(2\omega_4,0)\\[5pt]\displaybreak[0]
I(\omega_2,3)_\Gamma&=V(\omega_2,3)\oplus \left(\begin{array}{c}V(\omega_4,2)\\ \oplus\\ V(2\omega_2,2)\\ \oplus \\
                                     V(\omega_1+\omega_3)
                                    \end{array}\right)
\oplus 2 V(\omega_2+\omega_4,1)\oplus V(2\omega_4,0)\\\displaybreak[0]
I(\omega_4,2)_\Gamma&=V(\omega_4,2)\oplus V(\omega_2+\omega_4,1)\oplus V(2\omega_4,0)\\\displaybreak[0]
I(2\omega_2,2)_\Gamma&=V(2\omega_2,2)\oplus V(\omega_2+\omega_4,1)\oplus V(2\omega_4,0)\\\displaybreak[0]
I(\omega_1+\omega_3,2)_\Gamma&=V(\omega_1+\omega_3,2)\oplus V(\omega_2+\omega_4,1)\\\displaybreak[0]
I(\omega_2+\omega_4,1)_\Gamma&=V(\omega_2+\omega_4,1)\oplus V(2\omega_4,0)\\\displaybreak[0]
I(2\omega_4,0)_\Gamma&=V(2\omega_4,0).
\end{align*}

\subsection{}\label{app.2} We now establish the following result which was used in the proof of Proposition \ref{QR}. We use
the notation of the proof freely.
\begin{lem} Let $\lambda,\mu\in\{2\omega_2, \omega_4,
\omega_1+\omega_3\}$.  The elements $\bop\circ(1\otimes\bar
g_\mu)\circ \bar f_\mu$ and $\bop\circ(1\otimes\bar g_\lambda)\circ
\bar f_\lambda$ are linearly independent elements of $\Hom_{\lie
g}(V(\omega_2+\omega_4), S^2(\lie g)\otimes \lie g)$.
\end{lem}
\begin{pf}

Write~$x^-_i$ for~$x^-_{\alpha_i}$. Let~$w_{\omega_4}=\bar
g_{\omega_4}(v_{\omega_4})$. Since~$\lie g\cong_{\lie g}
V(\omega_2)$, we have
\begin{multline*}
   w_{\omega_4} = v_{\omega _2}\otimes x^-_2 x^-_1 x^-_3
    x^-_2 v_{\omega _2}+x^-_2 x^-_1 x^-_3 x^-_2
    v_{\omega _2}\otimes   v_{\omega _2}-x^-_2
    v_{\omega _2}\otimes x^-_1 x^-_3 x^-_2 v_{\omega
    _2}\\-x^-_1 x^-_3 x^-_2 v_{\omega
    _2}\otimes x^-_2 v_{\omega _2}+x^-_1
    x^-_2 v_{\omega _2}\otimes x^-_3 x^-_2
    v_{\omega _2}+x^-_3 x^-_2 v_{\omega
    _2}\otimes x^-_1 x^-_2 v_{\omega _2}\in S^2(\lie g).
\end{multline*}
Write $w_{\omega_4}=w_{\omega_4}^{(1)}\tensor w_{\omega_4}^{(2)}$ in Sweedler notation.
Then
$$
v_1=w_{\omega_4}\tensor v_{\omega_2},\quad v_2=v_{\omega_2}\tensor w_{\omega_4},\quad v_3=w_{\omega_4}^{(1)}\tensor v_{\omega_2}\tensor
w_{\omega_4}^{(2)}
$$
is a basis of $U=\{ v\in {\lie g^{\tensor 3}}_{\omega_2+\omega_4}\,:\, \lie n^+ v=0\}$.
Let $s_1=v_1$, $s_2=v_2+v_3$, $s_3=v_2-v_3$. Then~$s_1$, $s_2$ form a basis of $U\cap S^2(\lie g)\tensor\lie g$,
while $s_3$ spans $U\cap \bigwedge^2 \lie g\tensor\lie g$.

We have
\begin{multline*}
\bar f_{2\omega_2}(v_{\omega_2+\omega_4}) = -  v_{\omega _2}\otimes x^-_1 x^-_2 x^-_3
    x^-_2 v_{2 \omega _2}+2   v_{\omega
    _2}\otimes x^-_2 x^-_1 x^-_3 x^-_2 v_{2 \omega
    _2}-\\3 x^-_2 v_{\omega _2}\otimes x^-_1 x^-_3
    x^-_2 v_{2 \omega _2}+3 x^-_1 x^-_2 v_{\omega
    _2}\otimes x^-_3 x^-_2 v_{2 \omega _2}+3
    x^-_3 x^-_2 v_{\omega _2}\otimes x^-_1
    x^-_2 v_{2 \omega _2}\\-3 x^-_1 x^-_3 x^-_2
    v_{\omega _2}\otimes x^-_2 v_{2 \omega _2}+6
    x^-_2 x^-_1 x^-_3 x^-_2 v_{\omega _2}\otimes
    v_{2 \omega _2}
\end{multline*}
and
\begin{multline*}
 \bar f_{\omega_1+\omega_3}(v_{\omega_2+\omega_4})= v_{\omega _2}\otimes x^-_1 x^-_2 x^-_3
    v_{\omega _1+\omega _3}-2   v_{\omega _2}\otimes
    x^-_2 x^-_1 x^-_3 v_{\omega _1+\omega _3}\\+
    v_{\omega _2}\otimes x^-_3 x^-_2 x^-_1 v_{\omega
    _1+\omega _3}+2 x^-_2 v_{\omega _2}\otimes
    x^-_1 x^-_3 v_{\omega _1+\omega _3}-2
    x^-_1 x^-_2 v_{\omega _2}\otimes x^-_3 v_{\omega
    _1+\omega _3}-\\2 x^-_3 x^-_2 v_{\omega
    _2}\otimes x^-_1 v_{\omega _1+\omega _3}+2
    x^-_1 x^-_3 x^-_2 v_{\omega _2}\otimes
    v_{\omega _1+\omega _3}.
\end{multline*}
Furthermore,
$$
\bar g_{\omega_2}(v_{\omega_1+\omega_3})=
  v_{\omega _2}\otimes x^-_2 v_{\omega
    _2}-x^-_2 v_{\omega _2}\otimes   v_{\omega
    _2}.
$$
Finally, $\bar g_{2\omega_2}(v_{2\omega_2})=v_{\omega_2}\tensor v_{\omega_2}$.

The composite map
$$
V(\omega_2+\omega_4)\stackrel{\bar f_{2\omega_2}}\longrightarrow \lie g\tensor V(2\omega_2)
\stackrel{1\tensor \bar g_{2\omega_2}}{\longrightarrow} \lie g^{\tensor 3}
$$
sends $v_{\omega_2+\omega_4}$ to
\begin{multline*}
  m_1 = -  v_{\omega _2}\otimes x^-_2 v_{\omega
    _2}\otimes x^-_1 x^-_3 x^-_2 v_{\omega
    _2}-  v_{\omega _2}\otimes
    x^-_1 x^-_2 v_{\omega _2}\otimes
    x^-_3 x^-_2 v_{\omega _2}\\-  v_{\omega
    _2}\otimes x^-_3 x^-_2 v_{\omega
    _2}\otimes x^-_1 x^-_2 v_{\omega
    _2} +2   v_{\omega _2}\otimes
    x^-_2 v_{\omega _2}\otimes x^-_1 x^-_3 x^-_2
    v_{\omega _2}+  v_{\omega _2}\otimes
    x^-_2 x^-_1 x^-_3 x^-_2 v_{\omega _2}\otimes
    v_{\omega _2}\\+  v_{\omega _2}\otimes v_{\omega
    _2}\otimes x^-_2 x^-_1 x^-_3 x^-_2 v_{\omega
    _2}-3 x^-_2 v_{\omega _2}\otimes
    x^-_1 x^-_3 x^-_2 v_{\omega _2}\otimes v_{\omega
    _2}+x^-_2 v_{\omega _2}\otimes v_{\omega
    _2}\otimes x^-_1 x^-_3 x^-_2 v_{\omega
    _2}\\+3 x^-_1 x^-_2 v_{\omega
    _2}\otimes x^-_3 x^-_2 v_{\omega
    _2}\otimes v_{\omega _2}+x^-_1 x^-_2
    v_{\omega _2}\otimes v_{\omega _2}\otimes x^-_3
    x^-_2 v_{\omega _2}+3 x^-_3
    x^-_2 v_{\omega _2}\otimes x^-_1 x^-_2
    v_{\omega _2}\otimes v_{\omega _2}\\+x^-_3
    x^-_2 v_{\omega _2}\otimes v_{\omega _2}\otimes
    x^-_1 x^-_2 v_{\omega _2}-3
    x^-_1 x^-_3 x^-_2 v_{\omega _2}\otimes
    x^-_2 v_{\omega _2}\otimes v_{\omega
    _2}+x^-_1 x^-_3 x^-_2 v_{\omega _2}\otimes
    v_{\omega _2}\otimes x^-_2 v_{\omega _2}+\\6
    x^-_2 x^-_1 x^-_3 x^-_2 v_{\omega _2}\otimes
    v_{\omega _2}\otimes v_{\omega _2}
\end{multline*}
Furthermore, the image of~$v_{\omega_2+\omega_4}$ under the composite map
$$
V(\omega_2+\omega_4)\stackrel{\bar f_{\omega_1+\omega_3}}{\longrightarrow} \lie g\tensor V(\omega_1+\omega_3)
\stackrel{1\tensor\bar g_{\omega_1+\omega_3}}{\longrightarrow} \lie g^{\tensor 3}
$$
is
\begin{multline*}
m_2 = -2   v_{\omega _2}\otimes   v_{\omega
    _2}\otimes x^-_2 x^-_1 x^-_3 x^-_2 v_{\omega
    _2}+2   v_{\omega _2}\otimes
    x^-_2 x^-_1 x^-_3 x^-_2 v_{\omega _2}\otimes
      v_{\omega _2}+2 x^-_2 v_{\omega
    _2}\otimes   v_{\omega _2}\otimes
    x^-_1 x^-_3 x^-_2 v_{\omega _2}\\-2 x^-_2
    v_{\omega _2}\otimes x^-_1 x^-_3 x^-_2 v_{\omega
    _2}\otimes   v_{\omega _2}-2 x^-_1
    x^-_2 v_{\omega _2}\otimes   v_{\omega
    _2}\otimes x^-_3 x^-_2 v_{\omega _2}+2
    x^-_1 x^-_2 v_{\omega _2}\otimes
    x^-_3 x^-_2 v_{\omega _2}\otimes
    v_{\omega _2}\\-2 x^-_3 x^-_2 v_{\omega
    _2}\otimes   v_{\omega _2}\otimes
    x^-_1 x^-_2 v_{\omega _2}+2 x^-_3
    x^-_2 v_{\omega _2}\otimes x^-_1 x^-_2
    v_{\omega _2}\otimes   v_{\omega _2}+2
    x^-_1 x^-_3 x^-_2 v_{\omega _2}\otimes
      v_{\omega _2}\otimes x^-_2 v_{\omega
    _2}\\-2 x^-_1 x^-_3 x^-_2 v_{\omega
    _2}\otimes   x^-_2 v_{\omega _2}\otimes
      v_{\omega _2}
\end{multline*}
Finally, the composite map
$$
V(\omega_2+\omega_4)\stackrel{\bar f_{\omega_4}}{\longrightarrow} \lie g\tensor V(\omega_4)\stackrel{1\tensor \bar g_{\omega_4}}{\longrightarrow}
\lie g^{\tensor 3}
$$
maps $v_{\omega_2+\omega_4}$ to $m_3 = v_{\omega_2} \tensor w_{\omega_4}$.
In particular, this implies that all composite maps are non-zero. Furthermore,
it can be shown that
\begin{align*}
m_1&=3 s_1+s_2-2s_3\\
m_2&=2 s_1-s_2+s_3\\
2m_3&=s_2+s_3.
\end{align*}
In particular, $m_1$, $m_2$ and~$m_3$ are linearly independent.
Suppose that~$x_1 m_1+x_2 m_2+x_3 m_3\not=0$ has zero projection onto~$S^2(\lie g)\tensor \lie g$.
Using the above equations we conclude that~$x_i\not=0$, $i=1,2,3$.
\end{pf}

\end{document}